\documentclass[10pt]{article}

\usepackage[a4paper]{geometry} 

\usepackage{a4wide}
\usepackage{amssymb}
\usepackage{amsfonts}
\usepackage{amsmath}
\usepackage{graphicx}
\usepackage{tikz}
\usepackage{float}
\input xy
\xyoption{arrow} 
\xyoption{matrix}

\date{}

\newtheorem{proposition}{Proposition}[section]
\newtheorem{theorem}[proposition]{Theorem}
\newtheorem{lemma}[proposition]{Lemma}

\newtheorem{corollary}[proposition]{Corollary}

\def\GK{{\rm  GK}\,}
\def\Kdim{{\rm K.dim }\,}

\def\der{\partial }

\def\nFM0{{\nu }_{F,M_0}}
\def\nFN0{{\nu }_{F,N_0}}
\def\nGN0{{\nu }_{G,N_0}}

\def\N0{ {\bf N}_0 }

\def\t{\otimes}

\def\v{\varphi}
\def\ra{\rightarrow}

\def\Xpm{X^{\pm }}

\def\s{\sigma}
\def\Z{\mathbb{Z}}

\def\l1{{\lambda}_1}

\def\a{\alpha}
\def\a0{ {\alpha }_0}
\def\a1{ {\alpha }_1}

\def\l{\lambda}

\def\nFGM0{{\nu }_{F,G,M_0}}

\def\nFN0{{\nu}_{F,N_0}}

\def\sm{{\sigma}^m}

\def\sm1{{\sigma}^{-1}}

\def\smtp1{{\sigma}^{-t+1}}

\def\S1{S^{-1}}

\def\Xpm1{X^{\pm 1}_1}

\def\sPM1{{\sigma }^{\pm 1}}
\def\sMP1{{\sigma }^{\mp 1 }}

\def\d{\delta}

\def\di{{\rm d.ind}}

\def\L{\Lambda}

\def\G{\Gamma}

\def\OO{{\cal O}}
\def\CA{{\cal A}}

\def\CD{{\cal D}}

\def\Ytm1{Y^{t-1}}
\def\Yim1{Y^{i-1}}

\def\CG{{\cal G}}

\def\supp{{\rm supp }}

\def\Der{{\rm Der }}
\def\ad{{\rm ad }}
\def\dim{{\rm dim }}

\def\SL2Z{ {\rm SL}_2({\bf Z}) }

\def\CR{ {\cal R}}

\def\Gp1{ G^{1 , 1 } }
\def\P11{ P^{-1 , 1 } }
\def\Pp1{ P^{1 , 1 } }

\def\Supp{{\rm Supp}}

\def\nCLsr{{}^\nu\kern-2pt {\cal L}^{\sigma , \rho  }}
\def\nP{{}^\nu \kern-2pt P}
\def\nL{{}^\nu\kern-2pt L}
\def\nLL{{}^\nu\kern-2pt \Lambda}
\def\nPsr{{}^\nu\kern-2pt P^{\sigma , \rho  }}
\def\nLsr{{}^\nu\kern-2pt L^{\sigma , \rho  }}
\def\nuCL{{}^\nu\kern-2pt  {\cal L}}
\def\nCLsr{{}^\nu\kern-2pt {\cal L}^{\sigma , \rho  }}
\def\nCL1m{{}^\nu\kern-2pt {\cal L}^{-1 , 1  }}
\def\x1nu{x^\frac{1}{\nu}}
\def\xm1nu{x^{-\frac{1}{\nu}}}

\def\CR{ {\cal R}}

\def\ra{\rightarrow }

\def\CB{{\cal B}}

\def\nAM0{{\nu }_{{\cal A},M_0}}
\def\nAN0{{\nu }_{{\cal A},N_0}}

\def\Kdim{ {\rm Kdim } }

\def\End{ {\rm End }}
\def\Der{ {\rm Der }}

\def\CR{ {\cal R }}

\def\ad{ {\rm ad }}

\def\ga{\mathfrak{a}}

\def\gm{\mathfrak{m}}
\def\gp{\mathfrak{p}}
\def\gq{\mathfrak{q}}

\def\SL{{\rm SL}}

\def\Ext{{\rm Ext}}

\def\di!{\frac{\der^i}{i!}}
\def\dik!{\frac{\der^k_i}{k!}}

\def\gl{\mathfrak{l}}

\def\Max{{\rm Max}}

\def\gldim{{\rm gldim}}

\def\N{\mathbb{N}}

\def\0{\overline{0}}
\def\1{\overline{1}}

\def\Ln1{\L_{n,\overline{1}}}

\def\a1{a_{\overline{1}}}

\def\S{\Sigma}

\def\vn1{\overrightarrow{n-1}}

\def\gl{{\rm gl}}
\def\sl{{\rm sl}}

\def\mA{\mathbb{A}}

\def\mJ{\mathbb{J}}
\def\mI{\mathbb{I}}

\def\K1{{\rm K}_1}

\def\mA{\mathbb{A}}
\def\mB{\mathbb{B}}

\def\Supp{{\rm Supp}}

\def\hmI1{\widehat{\mI_1}}
\def\tmI1{\widetilde{\mI_1}}
\def\tmJ1{\widetilde{\mJ_1}}
\def\hB1{\widehat{B_1}}
\def\hCB1{\widehat{\CB_1}}

\def\ga{\mathfrak{a}}

\def\RSm1{R\langle S^{-1}\rangle}

\def\sl2{\mathfrak{sl}_2}
\def\gl2{\mathfrak{gl}_2}

\def\sl2{\mathfrak{sl}_2}

\def\CO{\mathcal{O}}

\def\Am{A(m)}
\def\CDAm{{\cal D} (A(m))}

\begin{document}

\author{V. V. \  Bavula  and K. Hakami 
}

\title{Rings of differential operators on singular generalized multi-cusp algebras}

\maketitle

\begin{abstract}

The aim of the paper is to study the ring of differential operators $\CDAm$ on the generalized multi-cusp algebra $A(m)$ where $m\in \N^n$ (of Krull dimension $n$). The algebra $\Am$ is singular apart from the single case when $m=(1, \ldots , 1)$. In this case, the algebra $\Am$ is a polynomial algebra in $n$ variables. So, the $n$'th Weyl algebra $A_n=\CD (A(1, \ldots , 1))$ is a member of the family of algebras $\CDAm$. 
We prove that the algebra $\CDAm$  is a central, simple, $\Z^n$-graded,  finitely generated Noetherian domain of Gelfand-Kirillov dimension $2n$. 
 Explicit finite sets of generators and defining relations is given for the algebra $\CDAm$. We prove that the Krull dimension and the global dimension of the algebra $\CDAm$ is $n$. An analogue of the Inequality of Bernstein is proven. In the case when $n=1$, simple $\CD (\Am)$-modules are classified. \\

 {\em Key Words:  the ring of differential operators, the generalized  multi-cusp algebra,  the generalized Weyl algebra, the global dimension, the Krull dimension, the Gelfand-Kirillov dimension, simple module, the projective dimension,  orbit,  the projective resolution.
 }\\

 {\em Mathematics subject classification 2010: 16S32, 16E10, 16E05, 16D60,  16P40,  16S35, 16S15, 16P50, 16P90, 16D30.}

$${\bf Contents}$$
\begin{enumerate}
\item Introduction.
\item Generators and defining relations of the algebra $\CD (A)$. 
\item Classification of simple $\CD (A)$-modules. 
 \item The global dimension of the algebra  $\CD (A)$.
 \item The algebras  $\CDAm$ where $m\in \N^n$.
\end{enumerate}
\end{abstract}


\section{Introduction}

The following notation will remain  fixed throughout the paper (if
it is not stated otherwise): $K$ is a field of  characteristic zero 
(not necessarily algebraically closed), module means a left
module, $P_n=K[x_1, \ldots , x_n]$ is a polynomial algebra over
$K$, $\der_1:=\frac{\der}{\der x_1}, \ldots ,
\der_n:=\frac{\der}{\der x_n}\in \Der_K(P_n)$,
$A_n=K\langle x_1, \ldots , x_n, \der_1, \ldots , \der_n \rangle \subseteq \End_K(P_n)$ is the $n$'th {\em Weyl algebra} ($A_n=\CD (P_n)$ is the algebra of differential operators on the polynomial algebra $P_n$),  
 $\N =\{ 0, 1, \ldots \}$ is the set of natural numbers and  $\N_{\geq s}=\{ i\in \N\, | \, i\geq s\}$. In the  case $n=1$, we usually drop the subscript `1'. So, $P=K[x]$ is a polynomial algebra in a variable $x$, 
$A_1=K \langle x , \der \; | \;  \der x - x \der = 1 \rangle $ 
is the {\em Weyl algebra}, i.e., $A_1=\CD (P)$ 
is the ring of differential operators on the polynomial algebra $P$.

The algebra of regular functions on the cusp $y^2=x^3$ is isomorphic to the subalgebra $A(2)=K + \sum_{i \geq 2} K x^{i}$ of the polynomial algebra $P=K[x]$. For each $m \geq 1$, 
$A=A(m)=K + \sum_{i \geq m} K x^i$ is a subalgebra of $P$ which is caleld the {\bf generalized cusp algebra}. Clearly, $A(1)=K[x]$ is a polynomial algebra and $A(2)$ is the cusp algebra. \\

 {\bf Definition.}  Let $m=(m_1, \ldots ,m_n)\in \N^n$, the subalgebra  of the polynomial algebra $P_n=K[x_1, \ldots , x_n]$, 
 $$\Am= \bigotimes_{i=1}^n A(m_i), \;\; {\rm where}\;\;  A(m_i) = K + \sum_{j \geq m_i} K x_i^j\subseteq K[x_i],$$ is called the {\bf generalized multi-cusp algebra} of rank $n$ (GMCA, for short).\\

Clearly, if $m=(1, \ldots , 1)$ then $\CD(A(m)) = A_n$ is the $n$'th Weyl algebra. 
 If if $m=(2, \ldots , 2)$ then $A(m)\simeq A(2)^{\t n}$ is the algebra of regular functions on the direct product of $n$ copies of the cusp.
 
 The aim of the paper is to study algebraic properties of the algebra $\CDAm$ 
of differential operators of the generalized multi-cusp algebra $\Am$ of rank $n$.  We are mostly interested in the case when $m=(m_1, \ldots , m_n)\in \N_{\geq 2}^n$ since for an arbitrary $m$ the algebra $\Am$ is isomorphic to the tensor product $P_s\t A(m')$ where $m'\in \N_{\geq 2}^{n-s}$ and $\CDAm\simeq A_s\t \CD (A(m'))$. \\

{\bf Generators and defining relations for the algebra $\CDAm$.} Theorem \ref{GenRel30Jan20} gives an explicit finite sets of generators and defining relations of the algebra $\CDAm$.

\begin{theorem}\label{GenRel30Jan20}
Let $m=(m_1, \ldots ,m_n)\in \N^n$.  Then 
\begin{enumerate}
\item $\CDAm\simeq \bigoplus_{i=1}^n\CD (A(m_i))$. 
\item For each $i=1, \ldots , n$, let $\CG_i$ and $\CR_i$ be the set of generators and defining relations of the algebra $\CD (A(m_i))$ as in  Theorem \ref{AB18Nov18}.(4). Then the algebra $\CDAm$ is generated by the finite set of elements $\CG = \cup_{i=1}^n\CG_i$ that satisfy the defining relations $\CR_1, \ldots , \CR_n$ and $g_ig_j=g_jg_i$ for all $g_i\in \CG_i$ and $g_j\in \CG_j$ for all $i\neq j$.  
\end{enumerate}
\end{theorem}

 A $K$-algebra $R$ is called {\em central} if its centre $Z(R)$ is equal to the field $K$. Theorem \ref{GenProp30Jan20} is about general properties of the algebra $\CDAm$. 
\begin{theorem}\label{GenProp30Jan20}
Let $m=(m_1, \ldots ,m_n)\in \N^n$.  Then The algebra $\CDAm$ is a central, simple, $\Z^n$-graded,  finitely generated Noetherian domain of Gelfand-Kirillov dimension $2n$. 
\end{theorem}

{\bf An analogue of the Inequality of Bernstein for the algebras  $\CD (\Am )$.} The starting point of the $\CD$-module theory is the {\em Inequality of Bernstein}: {\em For  all nonzero finitely generated $A_n$-modules} $M$, $\GK (M)\geq n$. 

\begin{theorem}\label{BerIn30JAn20}
Let $m=(m_1, \ldots ,m_n)\in \N^n$. For  all nonzero finitely generated $\CDAm$-modules $M$, $\GK (M)\geq n$. 
\end{theorem}

{\bf The Krull and global dimensions of the algebra $\CD (\Am )$.} The Krull dimension of the Weyl algebra $A_n$ is $n$, \cite{MR}.

\begin{theorem}\label{Kdim30JAn20}
Let $m=(m_1, \ldots ,m_n)\in \N^n$. The Krull  dimension of the algebra $\CDAm $ is $n$. 
\end{theorem}
 The global  dimension of the Weyl algebra $A_n$ is $n$, \cite{MR}.

\begin{theorem}\label{gldim30JAn20}
Let $m=(m_1, \ldots ,m_n)\in \N^n$. The global dimension of the algebra $\CDAm $ is $n$. 
\end{theorem}

{\bf Classification of simple $\CD (A)$-modules where $A=A(m)=K+\sum_{i\geq m}Kx^i$.}


The set $\widehat{\CD (A)}$ of isomorphism classes of simple $\CD (A)$-modules is a disjoint union of two subsets: the set of $D$-torsion and the set of $D$-torsion free simple $\CD (A)$-modules where $D=K[h]$ and $h=\der x$. The sets of simple $D$-torsion and  $D$-torsion free  $\CD (A)$-modules are classified in Theorem \ref{DTor27Jan20} and Theorem \ref{DTFr27Jan20}, respectively.


\section{Generators and defining relations of the algebra $\CD (A)$}\label{PRLMN}

The aim of this section is to find generators and defining relatyions of the algebra $\CD (A)$ of differential operators on the algebra $A=\Am$ (Theorem \ref{AB17Nov18}). It is proven that the algebra $\CD (A)$ is a central simple Noetherian domain of Gelfand-Kirillov dimension 2 (Theorem \ref{AB18Nov18}.(1)).
The Krull dimension of the algebra $\CD (A)$ is 1 (Theorem \ref {E29Jan20}).
 Furthermore, for all  nonzero left ideals $I$  of the algebra $\CD (A)$,  the $\CD (A)$-module $\CD (A)/I$ has finite length (Theorem  \ref{C29Jan20}). We introduce two generalized Weyl algebras $\mA$ and $\mB$ such that $\mA\subset \CD (A)\subset \mB = T^{-1}\mA \simeq  T^{-1}\CD (A)$. 
 The properties of the algebra $\CD (A)$ is a mixture of properties of the algebras $\mA$ and $\mB$. \\

{\bf Generalized Weyl algebras $D(\s , a)$ of rank 1,  \cite{Bav-GWA-FA-91}-\cite{Bav-Bav-Ideals-II-93}}. Let $D$ be a ring, $\s $ be a ring automorphism of $D$, 
 $a$ is a {\em central} element of $D$. 
 The {\bf generalized Weyl algebra} of rank 1 (GWA, for short) 
 $D(\s , a)=D[X,Y; \s , a]$ is a ring 
 generated by the ring $D$  and two elements $X$ and $Y$ 
 that are subject to the defining relations:
\begin{equation}\label{clGWA}
 Xd=\s (d) X\;\;  {\rm and} \;\;  Yd=\s^{-1} (d)Y\;\; {\rm for \; all} \;\; d\in D,  \;\;
 YX=a \;\; {\rm and} \;\;   XY=\s (a).
\end{equation}
The ring $D$ is called the {\em base ring} of the GWA. 
The automorphism $\s$ and the element $a$ are 
called the {\em defining automorphism} and 
the {\em defining element} of the GWA, respectively.

The algebra $A=\bigoplus_{n\in {\bf \Z}}\, A_n$
is $\Z$-graded where $A_n=Dv_n$,
$v_n=X^n$ and  $v_{-n}=Y^n$ for $n<0$,  and $v_0=1.$
 It follows from the above relations that
$v_nv_m=(n,m)v_{n+m}=v_{n+m}\langle n,m\rangle $
for some $(n,m)\in D$. If $n>0$ and $m>0$ then
\begin{eqnarray*}
 n\geq m:& &  (n,-m)=\sigma^n(a)\cdots \sigma^{n-m+1}(a),\;\;  (-n,m)=\sigma^{-n+1}(a)
\cdots \sigma^{-n+m}(a),\\
n\leq m: & & (n,-m)=\sigma^{n}(a)\cdots \sigma(a),\,\,\,\;\;\; \;\;\; \;\;\; (-n,m)=\sigma^{-n+1}(a)\cdots a,
\end{eqnarray*}
in other cases $(n,m)=1$. Clearly, $\langle n,m\rangle =\s^{-n-m}((n,m))$.\\

{\bf Example.} The Weyl algebra $A_1=K[h][x, \der ; \s , a=h]$ is a GWA where $h=\der x$ and  $\s (h)=h-1$. \\

{\bf Generators and defining relations of the algebra $\CD (A)$.} 
 The set $S_{x} = \{ x^{i} \; | \;  i \geq 0 \}$ 
(resp., $S_{x^m} = \{ x^{i  m} \; | \;  i \geq 0 \}$) 
is a multiplicative set of $P$ 
(resp., $P$ and $A$).
Clearly, 
\begin{equation}\label{SxmPA}
K[x,x^{-1}] = S_{x}^{-1} P = S_{x^m}^{-1} P = S_{x^m}^{-1} A.
\end{equation}

The polynomial algebra $P$ 
is a left $A_1$-module which is isomorphic to 
the factor module $A_1 / A_1 \der$ 
where the action of $A_1$ is given by the rule: 
For all $p \in P$, $x \ast p = xp$ and $\der \ast p = p := \frac{dp}{dx}$.
The left $A_1$-module $P=\bigoplus_{i \geq 0} K x^{i}$ is 
a $\Z$-graded 
(even $\N$-graded) $A_1$-module  and 
$h \ast x^{i} = (i +1) x^{i} \;\; {\rm  for \; all} \;\;  i \geq0 \;\; {\rm  where} \;\; h= \der x $.

\begin{theorem}\label{AB17Nov18}
Let $K$ be a field of characteristic zero, 
$A=K + \sum_{i \geq m} K x^{i}$ \; $(m \geq 2)$ be 
a subalgebra of the  polynomial algebra $P=K[x]$. 
Then 
\begin{enumerate}
\item The ring of differential operators $\CD (A)$ on $A$ 
is a $\Z$-graded subalgebra 
$ \CD (A) = \bigoplus_{i \in \Z} \CD (A)_{[i]}$ 
of the $\Z$-graded algebra $A_{1,x}$ where $\CD (A)_{[i]} = D \d_i$ and 
\[
\d_i = \begin{cases}
            x^{i} & \text{if} \, \, \,  i \geq m,  \\
             (h-i-1) x^{i}  & \text{if} \, \, \,  i=1, \ldots , m-1, \text{and} \\ 
             1  & \text{if} \, \, \, i=0, \\
             \end{cases}
\]
\[ 
\d_{-i} = \begin{cases}
            (h+i-1) \cdot \prod_{j=m-i+1}^m (h-j) x^{-i} & \text{if} \, \, \,  i = 1, \ldots, m-1,   \\ 
             (h+i-1) \cdot \prod_{1\neq j=m-i + 1}^m (h-j) x^{-i}  & \text{if} \, \, \, i \geq m. \\
             \end{cases}
\]
In particular, $\d_{-m}=(h+m-1) (h-2) \cdots (h-m) x^{-m}$, and for all $i\in \Z$, $\d_i=\v_ix^i$ where the polynomial $\v_i\in D=K[h]$ is  the coefficient of $x^i$ in the equalities above.  
\item For all $i,j \geq m$, $\d_{-i} \d_{-j} = \d_{-i-j}$ and $\d_i\d_j=\d_{i+j}$.
\item   $ \CD(A) = \bigoplus_{j\geq 0, m\leq i\leq 2m-1} D\d_{-i}\d_{-m}^j\oplus \bigoplus_{|i|<m}D\d_i\oplus \bigoplus_{j\geq 0, m\leq i\leq 2m-1} D\d_i\d_m^j $, and  $\d_{-i}\d_{-m}^j=\d_{-m}^j\d_{-i}$ and $\d_i \d_m^j=\d_m^j\d_i $ for all $j\geq 0$ and $ m\leq i\leq 2m-1$. 

\item The algebra $\CD(A)$ is  generated algebra by the elements 
$ \{ h , \d_i \; | \;  i=\pm 1, \pm 2, \ldots, \pm (2m-1) \}$ that satisfy the finite set of defining relations: For all $i,j=\pm 1, \ldots , \pm (2m-1)$,  $[h, \d_i]=i\d_i$  and 
$$
\d_i\d_j=\begin{cases}
\v_i\s^i (\v_j)\v_{i+j}^{-1} \d_{i+j}& \text{if } |i+j|<2m,\\
\v_i\s^i (\v_j)\v_{i+j-m}^{-1}\d_{i+j-m}\d_m & \text{if } 2m\leq i+j<3m,\\
\v_i\s^i (\v_j)\v_{i+j-2m}^{-1}\d_{i+j-2m}\d_m^2 & \text{if } 3m\leq i+j<4m,\\
\v_i\s^i (\v_j)\v_{i+j+m}^{-1}\d_{i+j+m}\d_{-m} & \text{if } -3m< i+j\leq -2m,\\
\v_i\s^i (\v_j)\v_{i+j+2m}^{-1}\d_{i+j+2m}\d_{-m}^2 & \text{if } -4m< i+j\leq -3m. \\
\end{cases}
$$
\end{enumerate}
\end{theorem}

{\it Proof}. 1. The set $S_x =\{ x^{i} \; | \;  i \geq 0 \}$ 
(resp., $S_{x^m} = \{ x^{mi} \; | \;   i \geq 0 \}$) 
is an Ore set of the Weyl algebra $A_1$ 
(resp., of $A_1$ and $\CD (A)$)
and 
\begin{equation}\label{ASAD}
A_{1,x} := S_{x}^{-1} A_1 = S_{x^m}^{-1} A_1 = S_{x^m}^{-1} \CD(P) \simeq  \CD (S_{x^m}^{-1} P) \stackrel{(\ref{SxmPA})} = \CD (S_{x^m}^{-1} A) \simeq  S_{x^m}^{-1} \CD (A).
\end{equation}
Recall that the Weyl algebra 
$A_1 = D [x, \der ; \s , a = h ]$ is GWA 
when $D = K[h]$, $\s (h) = h-1$ and $h:= \der x$. 
In particular, the Weyl algebra 
$A_1 = \bigoplus_{i \in \Z} D v_i$ is a $\Z$-graded algebra 
where $v_{0} := 1$, $v_{i} = x^{i}$ and $v_{-i} = \der^{i}$  for $i \geq 1$.

Since the elements of the Ore set $S_x$ 
are homogeneous elements of the algebra $A_1$, 
the localized algebra $A_{1,x} =S_x^{-1} A_1$ is also a $\Z$-graded algebra
$A_{1,x} = \bigoplus_{i \in \Z} D x^i$ (since $\der = h x^{-1}$). 
By (\ref{ASAD}),
$\CD (A) = \{ \d \in A_{1,x} \; | \; \d \ast A \subseteq A  \}$.

Since the algebra $A$ is a $\Z$-graded subalgebra of the polynomial algebra $P$, the algebra $\CD (A)$ is also $\Z$-graded,
\begin{equation}\label{ASAD1}
\CD (A) = \bigoplus_{i \in \Z} \CD (A)_{[i]}  \;\; {\rm  where} \;\;  \CD(A)_{[i]} = \CD(A) \cap D x^i = \{ \d \in D x^i  \; | \;  \d \ast A \subseteq A \}.
\end{equation}

Now, using the fact that $ h \ast x^i = (i+1) x^i$ 
for all $i \in \Z$, 
we obtain the explicit expressions for the graded components  
$\CD (A)_{[i]}$ as in the theorem.

2.  Statement 2 follows at once from the definition of the elements 
$\d_{-i}$ and $\d_i=x^i$ $(i \geq m)$  and the fact that 
$x^{-j} h = (h +j) x^{-j} $ for all $j \geq 0$.

3.  By statement 2, for all $i \geq 1$ and $j= 0,1, \ldots, m-1$, 
$\d_{-im -j} = \d_{-m}^{i} \d_{j}$ and $x^{im+j}=(x^m)^i \cdot x^j$.
Now, statement 3 follows from statement 1.

4. By statements 1 and 2, the relations in statement 4 hold. 
 Then the  relations of statement 4 are defining relations of the algebra $\CD (A)$ since they imply  the first equality in statement 3 where the direct sums are replaced by sums. $\Box$\\

{\bf The subalgebra $\CA$ of $\CD (A)$ which is a simple GWA.} 
 For an automorphism $\tau$ of an algebra $R$, $R^\tau := \{ r\in R\, | \, \tau (r)=r\}$ is the {\em algebra of $\tau$-constants/invariants}. The subalgebra $\CA$ of $\CD (A)$ (Proposition \ref{AB18Nov18}.(2)) plays a key role in proving that the algebra $\CD (A)$ is a simple Noetherian domain (Proposition \ref{AB18Nov18}.(1)).

\begin{theorem}\label{AB18Nov18}
Let $K$ be a field of characteristic zero, 
$A=K + \sum_{i \geq m} K x^{i} \;\; (m \geq 2)$ 
be a subalgebra of the  polynomial algebra $P=K [x]$.  Then 
\begin{enumerate}
\item The algebra $\CD (A)$ is a central simple  Noetherian domain. 
\item The subalgebra $\CA$ of $\CD (A)$ which is generated by the elements 
$h$, $X:= x^m$ and $Y:= \d_{-m}$ is a  GWA
 $\CA = D[X,Y;\s^m, a=(h+m-1) \cdot (h-2) (h-3) \cdots (h-m)]$ which  is a central simple  Noetherian domain
where $\s (h) = h-1$. 
\item The algebra $\CD(A)$ is a finitely generated  left and right $\CA$-module,  
$ \CD(A) = \sum_{|i| < 2m} \CA \d_i = \sum_{|i| < 2m} \d_{i} \CA $.
\end{enumerate}
\end{theorem}

{\it Proof}.  $2$.  The elements $h$, $X=x^m$ and $Y=\d_{-m}$ satisfy 
the defining relations for the GWA 
$D[X , Y ; \s^m , a]$. 
Then using the fact that the algebra $\CA$ is 
a homogeneous subalgebra of the $\Z$-graded algebra $\CD (A)$, 
we see that 
$\CA = \bigoplus_{i \geq 1} D Y^i \oplus D \oplus \bigoplus_{i \geq 1} D X^i$, and so $\CA = D [X , Y, \s^m , a]$ since $a=YX= (h+m-1) \cdot (h-2) (h-3) \cdots (h-m)x^{-m}x^m=(h+m-1) \cdot (h-2) (h-3) \cdots (h-m)$, by Theorem \ref{AB17Nov18}.(1). By  \cite[Corollary 3.2]{Bav-GWArep}, the GWA $\CA$  is a simple algebra since the difference of any two distinct roots of the polynomial 
 $a=(h+m-1) \cdot (h-2) (h-3) \cdots (h-m)$ is not divisible by $m$. 
 By  \cite[Proposition 1.3]{Bav-GWArep}, the GWA $\CA$  is a noetherian domain. Clearly, $Z(\CA ) = D^\s  = K$ since $\s (h)=h-1$ and the field $K$ has characteristic zero.

$3$.  By using the definition, 
the algebra $\CA$ is generated by the elements 
$h$, $x^m$ and $\d_{-m}$.
Now, statement 3 follows from Theorem \ref{AB17Nov18}.(2,3). 

$1$.   (i) {\em The algebra $\CD (A)$ is a Noetherian domain}: 
Since the polynomial algebra $D=K[h]$ is a Noetherian algebra, 
the GWA $\CA$ is also a  Noetherian domain \cite[Proposition 1.3]{Bav-GWArep}. 
The algebra $\CD(A)$ is a  finitely generated left and right $\CA$-module.
Hence, the algebra $\CD(A)$ is a Noetherian left and right $\CA$-module. 
Therefore, the algebra $\CD(A)$ is a Noetherian algebra. 

 (ii) {\em The algebra $\CD (A)$ is simple}: Let $I$ be a nonzero ideal of the algebra $\CD (A)$. Then $\ga := I\cap D\neq 0$ is a nonzero ideal of the algebra $D$ since the algebra $\CD (A)$ is a domain which is a direct sum $\CD (A)=\bigoplus_{i\in \Z} D\d_i$ of eigen-spaces of the inner derivation $\ad_h $ of the algebra $\CD (A)$ (Theorem \ref{AB17Nov18}.(1)). The subalgebra $\CA$ of $\CD (A)$ is a simple algebra (statement 2) that contains the algebra $D$. Hence, $0\neq \ga \subseteq I\cap \CA $ is a nonzero ideal of the algebra $\CA$, i.e. $1\in I\cap \CA\subseteq I$,  and so  $I=\CD (A)$. Therefore, the algebra $\CD (A)$ is a simple algebra. 

(iii) {\em The algebra $\CD (A)$ is central}: By (\ref{ASAD}), the algebra $\CD (A)$ is a central algebra:
$$K\subseteq Z(\CD (A))\subseteq Z(S_{x^m}^{-1}\CD (A))) =Z(A_{1,x}) =K.\;\; \Box$$

The set $S_{x^m}=\{ x^{im}\, | \, i\geq 0\}$ is a denominator set of the algebras $\CA$ and $A_1=\CD (P)$. The set $S_{x}=\{ x^{i}\, | \, i\geq 0\}$ is a denominator set of the Weyl  algebra $A_1=\CD (P)$. We have the following inclusions of algebras
\begin{equation}\label{AAA1xm}
 A_1\subset A_{1, x^m}=A_{1, x}=D[x,x^{-1}; \s], \;\; D=K[h],  \;\; \s (h) = h-1, 
\end{equation}
\begin{equation}\label{AAA1xm1}
 \CD (A )\subset \CD (A )_{ x^m}\simeq A_{1, x}=D[x,x^{-1}; \s],
\end{equation}
\begin{equation}\label{AAA1xm2}
\CA \subset \CA_{x^m}= D[x^m,x^{-m}; \s^m]\subset  A_{1, x^m}= A_{1, x}=D[x,x^{-1}; \s], 
\end{equation}
where the subscripts `$x^m$' and `$x$' denote the (left and right) localizations at the denominator sets $S_{x_m}$ and $S_x$, respectively. 
The rings  $D[x^{-1}, x; \s ]$ and $D[x^m,x^{-m}; \s^m]$ are  skew Laurent polynomial rings.

Recall that the Weyl algebra $\CD(P)=A_1$ is 
the GWA,  $A_1=D[x, \der ; \s, h] = \bigoplus_{i \in \Z} D v_i$,  
where $v_0=1$, $v_i = x^i$ and $v_{-i}= \der^{i}$ for all $i \geq 1$. 
Since $\der x = h$, 
we have that $ x^{-1} = h^{-1} \der$. 
Then, for all $i \geq 1$, 
\begin{equation}\label{xmiph}
x^{-i} = \prod_{k=0}^{i-1} (h+k)^{-1} \der^i.
\end{equation}
Now, for $i=1$, 
\begin{equation}\label{xmiph1}
\d_{-1} =\dfrac{h (h-m)}{h} \der = (h-m) \der.
\end{equation}

For $i=2, \ldots, m-1 $, 
\begin{equation}\label{xmiph2}
\d_{-i} =\dfrac{(h+i-1)\prod_{j=m-i+1}^{m} (h-j)}{\prod_{k=0}^{i-1} (h+k) } \der^i
 = \dfrac{\prod_{j=m-i+1}^{m} (h-j)}{\prod_{k=0}^{i-2} (h+k) } \der^i .
\end{equation}

For $i \geq m $, 
\begin{equation}\label{xmiph3}
\d_{-i} =\dfrac{(h+i-1) \prod_{1 \neq j = -i +m +1 }^{m} (h-j) }{ \prod_{k=0}^{i-1} (h+k) } \der^i 
= \dfrac{\prod_{j=2}^{m} (h-j)}{ \prod_{k=i-m}^{i-2} (h+k) } \der^i .
\end{equation}

\begin{corollary}\label{aAB17Nov18}
Let $A$ be as in Theorem \ref{AB17Nov18}.
\begin{enumerate}
\item $\CD (A) \not \subseteq \CD (P)$. 
\item Let $\d_{-i}$, $i \geq 1$ be as in Theorem \ref{AB17Nov18}. 
Then $\d_{-1} = \CD (P)$ and $\d_{-i} \not \in \CD (P)$ for $i \geq 2$. 
\end{enumerate} 
\end{corollary}

{\it Proof}. 
$1$. 
Statement 1 follows  from statement 2.

$2$. 
Statement 2 follows from (\ref{xmiph1}), (\ref{xmiph2}) and (\ref{xmiph3}). $ \Box$

In Theorem \ref{BB18Nov18}, the algebra 
$\CD(A) \cap \CD(P)$ is described and an explicit set of algebra generators is given for it. \\

{\bf The subalgebra $\CD (A)_{+}$ and $\CD (A)_{-} $ of $\CD (A)$}.
Let the algebra $A$ be as in Theorem \ref{AB17Nov18}. 
The algebra $\CD(A)$ contains two homogeneous subalgebras
$ \CD (A)_{+} := \bigoplus_{i \geq 0} D \d_i $ and 
$ \CD(A)_{-} := \bigoplus_{i \geq 0} D \d_{-i} $.

\begin{proposition}\label{CB18Nov18}
Let $A$ be as in Theorem \ref{AB17Nov18}. 
\begin{enumerate}
\item The algebras  $\CD(A)_{\pm} $ are finitely generated Noetherian algebras.
\item $\CD (A)_{+} \subseteq \CD (P)$ but $\CD (A)_{-} \not \subseteq \CD (P)$.
\item The algebra $ \CD (A)_{+} $ is a finitely generated, left and right module 
over its subalgebra $D [x^m ; \s^m]$ and 
the set $ \{ 1, \d_1, \ldots, \d_{2m-1} \}$ 
is a module generating set. 
\item The algebra $\CD (A)_{-}$ is a finitely generated, left and right module 
over its subalgebra $D [\d_{-m} ; \s^{-m}]$ and 
the set $ \{ 1, \d_{-1}, \ldots, \d_{-2m+1} \}$ is a module generating set.  
\end{enumerate} 
\end{proposition}

{\it Proof}. 
$2$. 
The inclusion $\CD (A)_{+} \subseteq \CD (P)$ is obvious. 
By Corollary \ref{aAB17Nov18}.(2), $\CD(A)_{-} \not \subseteq \CD(P)$.

$3$. 
Statement 3 follows at once from the explicit expressions for 
the elements $\d_i  \;\; (i \geq 0)$ and the fact that $\CD(A)_{+} = \bigoplus_{i \geq 0} D \d_i$.

$4$. 
Statement 4 follows from Theorem \ref{AB17Nov18}.(2,3).

$1$. 
The skew polynomial rings $D[x^m ; \s^m]$ and $D[\d_{-m} ; \s^{-m}]$ 
are Noetherian algebras (since $D$ is so).
Now, statement 1 follows at once from statement 4. $\Box$\\

{\bf The algebra $\CA_1$}. 
Recall that the algebras  $\CD (A)$ and $\CD(P)$ are 
homogeneous subalgebras of the $\Z$-graded algebra $A_{1,x}$. 
So, the intersection 
$\CA_1 : = \CD (A) \cap \CD(P) = \CD (A) \cap A_1$ 
is a   homogeneous subalgebra of the algebras  
$\CD(A)$, $A_1$ and $A_{1,x}$.
Clearly, 
$ \CA_1 = \{ \d \in \CD (P) \; | \; \d \ast A \subseteq A \}$.

\begin{theorem}\label{BB18Nov18}
Let the algebra $A$ be as in Theorem \ref{AB17Nov18}.
\begin{enumerate}
\item $ \CA_1 = \bigoplus_{i \in \Z} D w_i$ where $w_0 := 1$, $w_i = \d_i$, $w_{-i} = a_i \der^{i}$ for $i \geq 1$, and 
\[
a_i = \begin{cases}
             \prod_{j=m-i+1}^{m} (h-j) & \text{if} \, \, \,  i =1, \ldots, m-2,  \\
              \prod_{j=2}^{m} (h-j)  & \text{if} \, \, \,  i \geq m-1. \\ 
             \end{cases}
\]
\item The algebra $\CA_1$ is a finitely generated algebra and the set
 $\{ w_{-m} , w_{-1}=\d_{-1}, h, \d_1 , \ldots, \d_{m-1}, x^m  \}$
 is an algebra generating set, and $w_{-1}=\d_{-1}$.  *** find defining relations ***
 \item 
 			\begin{enumerate}
 					\item  For all $ i \geq m$, $w_{-i} w_{-1} = h w_{-i-1}$, $w_{-1} w_{-i} = (h - 1) w_{-1 -i}$ and $[w_{-i}, w_{-1}] = w_{-m -1}$. 
 					\item For $i = 1, \ldots, m-1$, $(w_{-1})^i = w_{-i}$. 
 					\item ${(w_{-1})}^{m} = h w_{-m}$.
 					\item For all $i \geq 1$, $[ \d_1 , x^i ] = i x^{i+1}$.
 					\item $w_{-1} \d_{1} = h (h-1)(h-m)$ and $\d_i w_{-1} = (h-1)(h-2)(h-m-1)= \s(h(h-1)(h-m))$.
 					\item For $i = 2, \ldots , m-1$, $w_{-1} \d_{i} = h (h-m) \d_{i-1}$ and $\d_{i} w_{-1} =  (h-i-1) (h-i-m) \d_{i-1}$.
 			\end{enumerate}
\end{enumerate}
\end{theorem}

{\it Proof}. 
$1$. 
Notice that $\CD(A)_{+} \subseteq A_1$, and so $\CD(A)_{+} \subseteq \CA_1$. 
Now statement 1 follows from the fact that the Weyl algebra 
$A_1 = D [ x, \der; \s, h] = \bigoplus_{i \geq 1} D \der^{i} \oplus D \oplus \bigoplus_{i \geq 1} D x^{i} $ is a GWA and
from (\ref{xmiph1})--(\ref{xmiph3}).

$3$. 
Straightforward.

$2$. 
By statement 1, the set 
$G = \{ h , w_i \; | \;  i \in \Z \backslash \{0\} \}$
is a generating set for the algebra $\CA_1$.
By the statements 3(a) and 3(b),
the elements $\{ w_i \; | \;  i \leq -m-1 \}$ and
$ \{ w_{-j} \; | \;  j = 2, \ldots , m-1 \}$ are redundant in $G$.
Similarly, by the statement 3(d), the elements 
$\{ w_{i} \; | \;  i \geq m+1 \}$ are also  redundant in $G$, 
and statement 2 follows. $\Box$\\

 {\bf The generalized Weyl algebras $\mA$ and $\mB$ such that $\mA \subset \CD (A)\subset \mB \subset T^{-1}\mA =T^{-1} \CD (A)= T^{-1} \mB $.} Let $\mA$ be the subalgebra of $\CD (A)$ which is generated by the elements $\d_{-1}$, $h$ and $\d_1$. By Theorem \ref{AB17Nov18}. (1), 
 $\d_{-1}= h(h-m)x^{-1}$ and $\d_1= (h-2)x$, and so the algebra 
\begin{equation}\label{ADdd}
\mA = D[\d_1, \d_{-1}; \s , h(h-1)(h-m)], \;\; D=K[h], \;\; \s (h) = h-1, 
\end{equation}
 is a GWA such that $\mA\subset \CA_1$ since $\d_{-1}, h , \d_1\in  \CA_1$ (Theorem \ref{BB18Nov18}.(2)). In particular, the algebra  $\mA = \bigoplus_{i>0}D\d_{-1}^i\oplus \bigoplus_{i\geq 0}D\d_1^i$ is a free left/right $D$-module where the set $\{ \d_{\pm 1}^i\, | \, i\geq 0\}$ is a free basis over $D$.

  The multiplicative submonoid $T= \langle h-i\, | \, i\in \Z\rangle$ of $D$ is a (left and right) denominator set of the algebras $\mA$, $
 \CA_1$, $\CD (A)$ and $A_1$ such that 
\begin{equation}\label{BBT}
T^{-1}\mA \simeq T^{-1}\CA_1  \simeq T^{-1} \CD (A) \simeq  T^{-1}A_1=:  \mB = T^{-1}D[x,x^{-1}; \s ], \; T^{-1}D=K[h, (h-i)^{-1}]_{i\in \Z}, \; \s (h) = h-1. 
\end{equation}
This follows from the explicit descriptions of the free bases over $D$ of the algebras $\mA$, $\CA_1$, $\CD (A)$ and $A_1=\bigoplus_{i>0}D\der^i\oplus \bigoplus_{i\geq 0}Dx^i$ (Theorem \ref{AB17Nov18}, Theorem \ref{BB18Nov18}). Notice that the algebra $\mB = T^{-1}D[x,x^{-1}; \s , 1] $ is a GWA where the ring $T^{-1}D$ is a Dedekind ring. 

 Similarly, the multiplicative set $D\backslash \{ 0\}$ is a (left and right) denominator set of the algebras $\mA$, $
 \CA_1$, $\CD (A)$ and $A_1$ such that 
\begin{equation}\label{BBT1}
D^{-1}\mA \simeq D^{-1}\CA_1  \simeq D^{-1} \CD (A) \simeq  D^{-1}A_1=:   B = K(h)[x,x^{-1}; \s ], \;\;\s (h) = h-1, 
\end{equation}
where $D^{-1}\mA$ denotes the localization  $(D\backslash \{ 0\})^{-1}\mA$ of the algebra $\mA$ at $D\backslash \{ 0\}$, and $K(h)$ is the field of rational functions in the variable $h$ over the field $K$. We have the following diagram of algebras where the vertical lines denote containments of the algebras:

\begin{figure}[h]
\centering
\begin{tikzpicture}
\node (A)  at (0,3) {$B= D^{-1}\mA = D^{-1}\CA_1  = D^{-1} \CD (A) =  D^{-1}A_1 = K(h)[x,x^{-1}; \s ]$};
\node (B)  at (0,2){$\mB= T^{-1}\mA = T^{-1}\CA_1  = T^{-1} \CD (A) =  T^{-1}A_1 = T^{-1}D[x,x^{-1}]$};
\node (C) at (-1,1) {$\CD(A)$};
\node (D)  at (1,1){$A_1$};
\node (E)   at (0, 0){ \qquad \qquad \qquad  $\CA_1=\CD (A) \cap A_1$};
\node (F)   at (-1, 0){$\CA$};
\node (G)   at (0,-1){$\mA$};
\draw [semithick, -] (A)--(B);
\draw [semithick, -] (B)--(C);
\draw [semithick, -] (B)--(D);
\draw [semithick, -] (C)--(F);
\draw [semithick, -] (E)--(G);
\draw [semithick, -] (C)--(E);
\draw [semithick, -] (D)--(E);
\end{tikzpicture} 
\end{figure}

We will see that the properties of the algebra $\CD (A)$ are a mixture of  properties of the GWAs $\mA$ and $\mB$. Theorem \ref{A24Jan20} and Theorem \ref{B24Jan20} are about some properties of the algebras $\mA$ and $\mB$. 

\begin{theorem}\label{A24Jan20}
\begin{enumerate}
\item The algebra $\mA$ is a finitely generated, central, non-simple Noetherian domain with $\GK (\mA ) = 2$. 
\item (\cite[Theorem 1.6]{Bav-1996THM}) $\gldim (\mA) =2$. 
\item (\cite[Theorem 2]{Bav-UkrMathJ-92}) All nonzero left ideals of the algebra $\mA$ are co-finite ($\dim_K(\mA / I)<\infty$). 
\item (\cite[Theorem 2]{Bav-UkrMathJ-92}) $\Kdim (\mA )=1$. 
\item (\cite[Theorem 4]{Bav-UkrMathJ-92}) In $\mA$ there are  only finitely many nonzero ideals. 
\item (\cite[Theorem 1]{Bav-UkrMathJ-92}) Up to isomorphism, there only tow simple finite dimensional $\mA$-modules: $L_1=\mA/\mA(\d_{-1}, h-1 ,\d_1)$, $\dim_K(L_1)=1$ and  
 $L_{m-1}=\mA/\mA(\d_{-1}^{m-1}, h-m ,\d_1)$, $\dim_K(L_{m-1})=m-1$.
 \item (\cite[Theorem 3.3]{Bav-GWArep}) The category of finite dimensional modules is not semisimple. 
 \item (\cite[Theorem 6]{Bav-GWArep}) For all simple $\mA$-modules $M$ and $N$,  the vector spaces $\Ext_{\mA}^i(M,N)$ and   ${\rm Tor}^{\mA}_i(M,N)$ are finite dimensional for all $i$.   
 \item (\cite[Theorem 4]{Bav-GWArep}) Let $M$ be a simple $\mA$-module and $q\in \mA \backslash K$, then the kernel and cokernel of the linear map $q_M :M\ra M$, $ m\mapsto qm$ are finite dimensional.
 
\end{enumerate}
\end{theorem}

\begin{theorem}\label{B24Jan20}
\begin{enumerate}
\item The algebra $\mB$ is a finitely generated, central, simple Noetherian domain with $\GK (\mB ) = 2$. 
\item (\cite[Theorem 1.6]{Bav-1996THM}) $\gldim (\mB) =1$. 
\item (\cite[Theorem 2]{Bav-UkrMathJ-92}) All nonzero left ideals of the algebra $\mB$ are co-finite ($\dim_K(\mB / I)<\infty$). 
\item (\cite[Theorem 2]{Bav-UkrMathJ-92}) $\Kdim (\mB )=1$. 
\item (\cite[Theorem 1, Theorem 5]{Bav-UkrMathJ-92}) All simple $\mB$-modules are infinite dimensional. 
\end{enumerate}
\end{theorem}

{\bf Every proper factor module of $\CD (A)$ has finite length and the Krull dimension of $\CD (A)$.}
 Recall that the algebra $\CD (A)$ is a finitely generated over its subalgebra $\CA$. Proposition shows that the subalgebra $\CA$ of $\CD (A)$ is large in the sense that it meets every nonzero left ideal of the algebra $\CD (A)$. 
 
\begin{proposition}\label{D29Jan20}
For all   nonzero left ideals $I$  of the algebra $\CD (A)$, $\CA \cap I\neq 0$.
\end{proposition}

{\it Proof}. The Gelfand-Kirillov dimensions of the domains $\CD (A)$ and $\CA$ is 2.  By Theorem \ref{AB18Nov18}, the algebra $\CD (A)$  is a finitely generated module over its subalgebra $\CA$. Hence, $2=\GK (\CA ) \leq \GK_\CA (\CD (A))\leq \GK (\CA ) =2$, and so $\GK_\CA (\CD (A))=2$. Then, by \cite[Proposition 8.3.5]{MR}, $$\GK_{\CA}(\CD (A)/I)<\GK_{\CA}(\CD (A))-1=2-1=1.$$
Hence, $\CA \cap I\neq 0$  since $\GK (\CA )=2>1=\GK_{\CA}(\CD (A)/I)$.  $\Box $

\begin{theorem}\label{C29Jan20}
For all   nonzero left ideals $I$  of the algebra $\CD (A)$,  the $\CD (A)$-module $\CD (A)/I$ has finite length.
\end{theorem}

{\it Proof}.  By \cite[Theorem 2.1]{Bav-GWArep}, for all   nonzero left ideals $I'$  of the algebra $\CA$,  the $\CA$-module $\CA /I'$ has finite length.  By Theorem \ref{AB18Nov18}, the algebra $\CD (A)$ is a finitely generated $\CA$-module. Now, the theorem follows from Proposition \ref{D29Jan20}.  $\Box $

\begin{theorem}\label{E29Jan20}
The Krull dimension of the algebra $\CD (A)$ is 1.
\end{theorem}

{\it Proof}. The theorem follows at once from Theorem \ref{C29Jan20}. $\Box$


\section{Classification of simple $\CD (A)$-modules}\label{SIMMOD}

The aim of this section is to classify simple $\CD (A )$-modules 
 where $A= K+\sum_{i\geq m} Kx^i$ (Theorem \ref{DTor27Jan20} and Theorem \ref{DTFr27Jan20}). They are partitioned in two (disjoint) sets: $D$-torsion and $D$-torsion free. 
 The simple $\CD (A)$-modules in each of the two sets are classified (Theorem \ref{DTor27Jan20} and Theorem \ref{DTFr27Jan20}). 

At the beginning of the section we recall a classification of simple modules over a generalized Weyl   $A=D(\s, a)=D[x,y;\s,a]$ where $D$ is  a (commutative) Dedekind domain with some extra condition on the automorphism that is satisfied for our GWAs. In all the papers we cite below these algebras are denoted by `$A$', we hope that this notation will not lead to confusion.  
 
For an algebra $A$, we denote by $\widehat{A}$ the set of isomorphism classes of simple $A$-modules.  For an $A$-module $M$,  we denote by $[M]$ its isomorphism class. If $P$ is an isomorphism invariant  property of simple modules  (e.g., `being weight') then $\widehat{A}\,(P)$ stands for the set of all isomorphism classes of simple $A$-modules that satisfy $P$.\\

{\bf Classification of simple $A$-modules where $A=D(\s, a)$ and $D$ is a Dedekind ring}.
Let $A=D(\s, a)=D[x,y;\s,a]$ be a GWA such that $D$ is a Dedekind ring, $a \neq 0$, and the automorphism $\s$ of $D$ satisfies the condition:\\

 (*)  $\s^i(\gp) \neq \gp$ {\em  for all $i \in\Z \setminus \{0\}$ and all maximal ideals $\gp$ of} $D$.\\
 
 \textbf{Example}. The Weyl algebra $A_1=K[h][x, \der ; \s ,h]$ is an example of the GWA $A$.\\

\textbf{Example}. $A=K[h] ( \s,a)$ where $\s (h)=h-1$ and $K$ is a field of characteristic zero. In particular, the algebras $\mA$ is  of this type, see (\ref{ADdd}). A classification of simple $K[h](\s,a)$-modules is given in \cite{Bav-UkrMathJ-92,Bav-GWArep}.\\

\textbf{Example}. The GWA $\mB$ is an example of the GWA $A$, see (\ref{BBT}).\\

The set $S:=D\setminus \{ 0\}$ is an Ore set of the domain $A$. So, a simple $A$-module $M$  is either $D$-{\em torsion} ($S^{-1}M=0$) or   $D$-{\em torsion free} ($S^{-1}M \neq 0$). In the second case, the $S^{-1}A$-module $S^{-1}M$ is simple. 

Let us recall a classification of simple $A$-modules for the algebra $A=D(\s,a)$, see \cite{Bav-GWArep,Bav-UkrMathJ-92,Bav-GWAKernel-1993} for details. Clearly,
\begin{equation}\label{AADT}
\widehat{A}= \widehat{A}\,(D\text{-}{\rm torsion}) \,\,\coprod\,\,\widehat{A}\,(D\text{-}{\rm torsion \; free}).
\end{equation}

\textbf{The set $\widehat{A}\,(D\text{-}{\rm torsion})=\widehat{A}\,({\rm weight})$}.
The group $\langle \s \rangle \simeq \Z$ acts freely on the set $\Max\,(D)$ of maximal ideals of the Dedekind ring $D$. For each maximal ideal $\gp$ of $D$, $\CO(\gp)=\{ \s^i(\gp)\,|\, i \in \Z \}$ is its orbit. We use the bijection $\Z \rightarrow \CO(\gp)$, $i \mapsto \s^i(\gp)$, to define the order $\leqslant$ on each orbit $\CO(\gp)$: $\s^i(\gp) \leqslant \s^j(\gp)$ iff $i \leqslant j$. A maximal ideal of $D$ is called \emph{marked} if it contains the element $a$. There are only finitely many marked ideals. An orbit $\CO$ is called \emph{degenerated} if it contains a marked ideal. Marked ideals, say $\gp_1 < \cdots <\gp_s$, of a degenerated orbit $\CO$ partition it into $s+1$ parts,
\begin{equation}\label{GGG}
 \G_1=(-\infty, \gp_1], \,\,\,\,  \G_2=(\gp_1, \gp_2], \ldots, \G_s=(\gp_{s-1}, \gp_s], \,\,\,\,\G_{s+1}=(\gp_s, \infty).
\end{equation}
Two ideals $\gp, \gq \in \Max\,(D)$ are called \emph{equivalent} $\gp \sim \gq$ if they belong either to a non-degenerated orbit or to some $\G_i$. We denote by $\Max\,(D)/\!\!\sim $ the set of equivalence classes in $\Max\,(D)$.

An $A$-module $V$ is called \emph{weight} if $V=\bigoplus_{\gp \in \Max\,(D)}V_{\gp}$ where $V_{\gp}=\{v \in V\, |\, \gp v =0 \} =$ the sum of all simple $D$-submodules of $V$ which are isomorphic to $D/\gp$. The set $\Supp\,(V)=\{ \gp \in \Max\,(D)\,|\, V_{\gp} \neq 0 \}$ is called the \emph{support} of $V$, elements of $\Supp\,(V)$ are called \emph{weights} and $V_{\gp}$ is called the \emph{component} of $V$ of weight $\gp$. Clearly, an $A$-module is weight iff it is a semisimple $D$-module. Clearly,
\begin{equation}\label{AD=AW}
\widehat{A}\,(D\text{-}{\rm torsion})=\widehat{A}\,({\rm weight}),
\end{equation}
i.e., a simple $A$-module is $D$-torsion iff it is weight.

\begin{theorem}\label{AA27Feb16} 
\cite{Bav-UkrMathJ-92,Bav-GWArep,Bav-GWAKernel-1993} 
{\sc (Classification of simple $D$-torsion/weight $A$-modules)}  \\
The map $\Max\,(D)/\!\!\!\sim  \,\, \rightarrow \widehat{A}\,(D\text{-}{\rm torsion})$, $\G \mapsto [L(\G)]$, is  a bijection with the inverse $[M] \mapsto \Supp\,(M)$ where
\begin{enumerate}
\item If $\G$ is  a non-degenerated orbit then $L(\G)=A/A\gp$ where $\gp \in \G$.
\item If $\G=(-\infty, \gp]$ then $L(\G)=A/A(\gp, x)$.
\item If $\G=(\s^{-n}(\gp), \gp]$ for some $n \geqslant 1$ then $L(\G)=A/A(y^n, \gp ,x)$. The $D$-length of $L(\G)$ is $n$.
\item If $\G=(\gp, \infty)$ then $L(\G)=A/A(\s(\gp), y)$.
\end{enumerate}
\end{theorem}

\textbf{The set $\widehat{A}\,(D\text{-}{\rm torsionfree})$}.
For elements $\alpha, \beta \in D$, we write $\alpha < \beta$ if $\gp < \gq$ for all $\gp, \gq \in \Max\,(D)$ such that $\CO(\gp)=\CO(\gq)$, $\alpha \in \gp$ and $\beta \in \gq$. (We write also $\alpha < \beta$ if there are no such ideals $\gp$ and $\gq$). Recall that the GWA $A=\bigoplus_{i \in \Z}A_i$ is a $\Z$-graded algebra where $A_i=Dv_i=v_i D$, $v_0=1$, $v_i=x^i$ and $v_{-i}=y^i$ for all $i\geq 1$.\\

\emph{Definition}, \cite{Bav-UkrMathJ-92,Bav-GWArep,Bav-GWAKernel-1993}.
An element $b=v_{-m}\beta_{-m}+v_{-m+1}\beta_{-m+1}+\cdots +\beta_0 \in A$ (where $m \geqslant 1$, all $\beta_i \in D$ and $\beta_{-m}, \beta_0 \neq 0$) is called a \emph{normal} element if $\beta_0 < \beta_{-m}$ and $\beta_0 < a$.\\

The set $S:=D\setminus \{ 0\}$ is an Ore set of the domain $A$. Let $k:=S^{-1}D$ be the field of fractions of $D$. The algebra $B:=S^{-1}A=k[x, x^{-1};\s]$ is a skew Laurent polynomial ring which is a (left and right) principle ideal domain. So, any simple $B$-module is of type $B/Bb$ for some irreducible element $b$ of $B$. Two simple $B$-modules are isomorphic, $B/Bb \simeq B/Bc$, iff the elements $b$ and $c$ are \emph{similar} (i.e., there exists an element $d \in B$ such that $1$ is the greatest common right divisor of $c$ and $d$, and $bd$ is a least common left multiple of $c$ and $d$).

\begin{theorem}\label{BA27Feb16} 
\cite{Bav-UkrMathJ-92,Bav-GWArep,Bav-GWAKernel-1993} {\sc (Classification of simple $D$-torsionfree $A$-modules)}\\
$\widehat{A}\,(D\text{-}{\rm torsionfree})=\{ [M_b:=A/A\cap Bb]\,|\,b$ is a normal irreducible element of $B\}$. The $A$-modules $M_b$ and $M_{b'}$ are isomorphic iff the elements $b$ and $b'$ are similar.
\end{theorem}

For all nonzero elements $\alpha, \beta \in D$, the $B$-modules $S^{-1}M_b$ and $S^{-1}M_{\beta b \alpha^{-1}}$ are isomorphic. If an element $b=v_{-m}\beta_{-m}+\cdots +\beta_0$ is irreducible in $B$ but not necessarily normal the next lemma shows  that there are explicit elements $\alpha$ and $\beta$ such that the element $\beta b \alpha^{-1}$ is normal and irreducible in $B$.

\begin{lemma}\label{a27Feb16} 
\cite[Lemma 13]{Bav-UkrMathJ-92} {\sc (Normalization procedure)} \\
Given an element $b= v_{-m} \beta_{-m}+\cdots+\beta_0 \in A$ where $m \geqslant 1$, all $\beta_i \in D$ and $\beta_{-m}, \beta_0 \neq 0$. Fix a natural number $s \in \N$ such that $\s^{-s} (\beta_0) < \beta_{-m}$, $\s^{-s}(\beta_0)< \beta_0$ and $\s^{-s}(\beta_0)<a$. Let $\alpha=\prod_{i=0}^s \s^{-i}(\beta_0)$ and $\beta=\prod_{i=1}^{s+m}\s^{-i}(\beta_0)$. Then the element $\beta b \alpha^{-1}$ is a normal element which is called a normalization of $b$ and denoted $b^{\rm norm}$ (we can always assume that $s$ is the least possible).
\end{lemma}
 
{\bf The algebra $A=K+\sum_{i\geq m} Kx^i$ is a simple weight $\CD (A)$-module.} Clearly, $A=\sum_{i\in E}Kx^i$  where $E:=\{ 0,  m, m+1, \ldots \}$.   Recall that  by the very definition of the algebra $\CD (A)$ of differential operators on $A$, the algebra $A$ is a left $\CD (A)$-module and the action of  an element $\d \in \CD (A)$ on an element $a\in A$ is denoted either by $\d *(a)$ or $\d (a)$. For all $i\in E$, $h*x^i=(i+1)x^i$. This implies that the $\CD (A)$-module $A$ is a weight $\CD (A)$-module with $\Supp (A)=\{ (h-i-1)\, | \, i\in E\}$. In particular the $\CD (A)$-module $A$ is $D$-torsion. 
If follows from the equalities 
\begin{equation}\label{Eldd}
\d_1= (h-2)x, \; \d_{-1}= h(h-m)x^{-1}, \; \d_{-1}\d_1= h(h-1)(h-m)\; {\rm and}\; \d_1\d_{-1}=  (h-1)(h-2)(h-m-1)
\end{equation}
that the maps
\begin{eqnarray*}
\d_1 &:& Kx^i\ra Kx^{i+1}, \;\;\;  p\mapsto \d_1*p, \; \;\;\; i\geq m, \\
 \d_{-1} &:& Kx^{i+1}\ra Kx^i, \;\;\; p\mapsto \d_{-1}*p, \; \; i\geq m,
\end{eqnarray*}
are bijections. Similarly, it follows from the equalities
$\d_{-m}\d_m =(h+m-1) (h-2) \cdots (h-m)$ and $\d_m \d_{-m}=(h-1) (h-m-2) \cdots (h-2m)$ that the maps
\begin{eqnarray*}
\d_m &:& K\ra Kx^m, \;\;  p\mapsto \d_m*p, \\
 \d_{-m} &:& Kx^m\ra K, \; \; p\mapsto \d_{-m}*p,
\end{eqnarray*}
are bijections. Therefore, the algebra $A$ is a simple weight $\CD (A)$-module  with $\Supp (A)=\{ (h-i-1)\, | \, i\in E\}$. \\

{\bf Classification of simple weight $\CD (A)$-modules with support that belongs to the orbit $\OO (h)$.} The ideal $(h) = Dh$ is a maximal ideal of the polynomial algebra $D=K[h]$ with $D/(h) =K$. Let $\OO (h)=\OO ((h))=\{ \s^i (h)=(h-i) \, | \, i\in \Z\}$ be its $\s$-orbit. We will see that (up to isomorphism) there are only two simple weight $\CD (A)$-modules with support in $\OO (h)$: the algebra $A$ and a `complementary' module $A'$ which we are going to define. Furthermore, $\supp (A') = \OO (h) \backslash \Supp (A)$.

The polynomial algebra $K [x]$ has the canonical structure of the left $A_1$-module. 
Namely, $K [x] \simeq A_1 / A_1 \der $; $ x \ast p = x p $ and 
$ \der \ast p = \dfrac{d p}{d x}$ for all $p \in P$.
The {\em Laurent polynomial algebra} $L = K [x , x^{-1}]=\bigoplus_{i\in \Z} Kx^i$,
which is the localization of the polynomial algebra $K[x]$ at 
$S_x = \{ x^{i} \; | \;  i \geq 0 \}$,
is a left $A_{1,x}$-module. 
By (\ref{ASAD1}), the Laurent polynomial algebra $L$ is a left module over the algebras $A_{1,x}\simeq S_{x^m}^{-1}\CD (A)$ and $\CD (A)$. One can easily verifies using Theorem \ref{AB17Nov18}, that the subalgebra $A$ is a $\CD (A)$-submodule of $L$. Consider the $\CD (A)$-module, 
\begin{equation}\label{Apdef}
A':=L/A=\bigoplus_{i\in E'}Kx^i, \;\; E':=\Z\backslash E=\{ \ldots , -2,-1, 1, 2, \ldots , m-1\}.
\end{equation}
By (\ref{Eldd}),  the maps
\begin{eqnarray*}
\d_1 &:& Kx^i\ra Kx^{i+1}, \;\;  p\mapsto \d_1*p, \; \;\;\;\; i\in E'\backslash \{ -1, m-1\}, \\
 \d_{-1} &:& Kx^{i+1}\ra Kx^i, \;\; p\mapsto \d_{-1}*p, \;\;\;  i\in E'\backslash \{ -1, m-1\},
\end{eqnarray*}
are bijections.
Since $\d_2= (h-3)x^2$ and $\d_{-2}= (h+1) (h-m+1) (h-m) x^{-2}$, we have that 
\begin{equation}\label{Eldd1}
 \d_{-2}\d_2= (h+1) (h-m+1) (h-m) (h-1)
\; {\rm and}\; \d_2\d_{-2}= (h-3) (h-1)(h-m-1)(h-m-2), 
\end{equation}
and so   the maps
\begin{eqnarray*}
\d_2 &:& Kx^{-1}\ra Kx, \;\;  p\mapsto \d_2*p, \\
 \d_{-2} &:& Kx\ra Kx^{-1}, \; \; p\mapsto \d_{-2}*p,
\end{eqnarray*}
are bijections.  Therefore, the $\CD (A)$-module   $A'$ is a simple weight $\CD (A)$-module  with $\Supp (A')=\OO (h) \backslash \Supp (A)= \{ (h-i-1)\, | \, i\in E'\}$. 

\begin{lemma}\label{a26Jan20}
The $\CD (A)$-modules $A$ and $A'$ are the only two (up to isomorphism) simple weight $\CD (A)$-modules with support in the orbit $\OO (h)$. 
\end{lemma}

{\it Proof}. Recall that the $\CD (A)$-modules $A$ and $A'$ are non-isomorphic   simple weight $\CD (A)$-modules with support in the orbit $\OO (h)$. Now, the lemma follows at once from the fact that every simple weight module is uniquely determined by its support and that $\OO (h)=  \Supp (A) \coprod \Supp (A')$. $\Box $

Let us collect properties of the $\CD (A)$-modules $A$ and $A'$ in the next two lemmas. 

\begin{lemma}\label{b26Jan20}
\begin{enumerate}
\item The algebra $A=\bigoplus_{i\in E}Kx^i$ is a simple weight $S_{x^m}$-torsion $\CD (A)$-module with  $\Supp (A)=\{ (h-i-1)\, | \, i\in E\}$  where  $E=\{ 0,  m, m+1, \ldots \}$, and $\End_{\CD (A)}(A)=K$. 
\item ${}_{\CD (A)}A\simeq \CD (A) / \CD (A) (h-1, \d_{-1})=\bigoplus_{i\in E} K\d_i\overline{1} $ where $\overline{1}:= 1+\CD (A) (h-1, \d_{-1})$. 
\end{enumerate}
\end{lemma}

{\it Proof}. 1. The weight spaces of the weight $\CD (A)$-module $A$  are 1-dimensional,  hence  $\End_{\CD (A)}(A)=K$. The rest of statement 1 have been proven above.

2. The $\CD (A)$-module $W=\CD (A) / \CD (A) (h-1)\simeq \bigoplus_{i\in \Z }K\d_i1^*$ is weight  with $\Supp (W)=\Z$ where $1^*=1+\CD (A) (h-1)$. The map $W\ra A$, $1^*\mapsto \overline{1}$ is a $\CD (A)$-module epimorphism. Hence, ${}_{\CD (A)}A\simeq \CD (A) / \CD (A) (h-1, \d_{-1})$, by Lemma \ref{a26Jan20}.   $\Box $

\begin{lemma}\label{c26Jan20}
\begin{enumerate}
\item The algebra $A'=\bigoplus_{i\in E'}Kx^i$ is a simple weight $S_{x^m}$-torsion $\CD (A)$-module with  $\Supp (A')=\OO (h)\backslash \Supp (A)= \{ (h-i-1)\, | \, i\in E'\}$ 
 where  $E'= \Z \backslash E =  \{ \ldots , -2,-1, 1, 2, \ldots , m-1\}$, and $\End_{\CD (A)}(A')=K$. 
\item ${}_{\CD (A)}A'\simeq \CD (A) / \CD (A) (h, \d_1)=\bigoplus_{i\in E'} K\d_i\overline{1}' $ where $\overline{1}':= 1+\CD (A) (h, \d_1)$. 
\end{enumerate}
\end{lemma}

{\it Proof}. 1. The weight spaces of the weight $\CD (A)$-module $A'$  are 1-dimensional,  hence  $\End_{\CD (A)}(A')=K$. The rest of statement 1 have been proven above.

2. The $\CD (A)$-module $W'=\CD (A) / \CD (A) h\simeq \bigoplus_{i\in \Z }K\d_i1^o$ is weight  with $\Supp (W)=\Z$ where $1^o=1+\CD (A) h$. The map $W'\ra A'$, $1^o\mapsto \overline{1}'$ is a $\CD (A)$-module epimorphism. Hence, ${}_{\CD (A)}A'\simeq \CD (A) / \CD (A) (h, \d_1)$, by Lemma \ref{a26Jan20}.   $\Box $\\

{\bf Classification of simple $D$-torsion $\CD (A)$-modules.}   Recall that $\mA\subset \CD (A) \subset \mB = T^{-1}\mA=T^{-1}\CD (A)$. 
 So, every $\mB$-module is automatically is an $\mA$-module and $\CD (A)$-module. The group $\langle \s \rangle $ acts on the set $\Max (D)$ of maximal ideal of the algebra $D=K[h]$. The field $K$ has characteristic zero and $\s (h) = h-1$. So, every orbit $\OO (\gp ) = \{ \s^i (\gp )\, | \, i\in \Z\}$ contains infinite number of elements where $\gp \in \Max (D)$. 
 We denote by $ \Max (D)/\langle \s \rangle $ is the set of all $\s$-orbits in  $ \Max (D)$.

 The algebra $\mB = T^{-1}D[x,x^{-1}; \s , 1]$ is a GWA where $ T^{-1}D$ is a Dedekind ring and the automorphism $\s$ satisfies the condition (*) above. 
 Notice that $\Max (T^{-1} D) = \{ T^{-1}\gp \, | \, \gp \in \Max (D)\backslash \OO (h)\} $ where $\OO (h)$ is the $\s$-orbit of the maximal ideal $(h)$ of the algebra $D$, and the map 
 $ \Max (D)\backslash \OO (h)\ra \Max (T^{-1} D)$, $\gp\mapsto T^{-1}\gp$ 
 is a bijection. 
 
 For each orbit $\OO\in \Max (D)/\langle \s \rangle \backslash \{ \OO (h)\}$, we fix its element, say $\gp_\OO$. So, $\OO (\gp_\OO )=\OO $.
 \begin{proposition}\label{A29Jan20}

\begin{enumerate}
\item $\widehat{\mB} (T^{-1}D-{\rm torsion})= \{ \mB / \mB \gp_\OO  \, | \, \gp_\OO \in \Max (D)/\langle \s \rangle \backslash \{ \OO (h)\} \}$. 
\item The restriction map  $\widehat{\mB} (T^{-1}D-{\rm torsion})\ra  \widehat{\CD (A)} (D-{\rm torsion})$, $M\ra {}_{\CD (A)}M$ is an injection. 
\end{enumerate}
\end{proposition}

{\it Proof}. 1. Statement 1 follows at once from Theorem \ref{AA27Feb16} and the fact that the defining element of the GWA $\mB$ is $1$, and so every orbit of the automorphism $\s$ in $\Max (T^{-1}D)$ is not degenerated.

2. Given $[M]\in \widehat{\mB} (T^{-1}D-{\rm torsion})$. By statement 1,  $$M= \mB / \mB \gp =\bigoplus_{i\in \Z} x^{-i}T^{-1}D/T^{-1}D\gp \simeq \bigoplus_{i\in \Z} x^{-i} D/D\gp $$ is a direct sum of non-isomorphic simple $D$-modules  for some $\gp = \gp_\OO \in \Max (D)\backslash \OO (h)$. By  Theorem \ref{AA27Feb16} in case of the GWA $\mA$, the weight $\mA$-module $M$ is simple, hence the $\CD (A)$-module $M$ is simple since $\mA\subset \CD (A)$. $\Box $
 
 In view of Proposition \ref{A29Jan20}.(2), we can write $\widehat{\mB} (T^{-1}D-{\rm torsion}) \subseteq  \widehat{\CD (A)} (D-{\rm torsion})$.

 \begin{theorem}\label{DTor27Jan20}
{\sc (Classification of simple $D$-torsion $\CD (A)$-modules)} 
\begin{enumerate}
\item $\widehat{\CD (A)} (D-{\rm torsion})=\{ A,  A'\} \coprod \widehat{\mB} (T^{-1}D-{\rm torsion})$. 
\item  $\widehat{\CD (A)} (D-{\rm torsion})=\{ A,  A'\}  \coprod \{ \mB / \mB \gp_\OO  \, | \, \gp_\OO  \in \Max (D)\backslash \OO (h))\}$, $\Supp (\mB / \mB \gp_\OO )= \OO $. 
\item For all $[M]\in \widehat{\CD (A)} (D-{\rm torsion})$, $\Supp (M)=\infty$ and $\dim_K(M) = \infty$. 
\end{enumerate}
\end{theorem}

{\it Proof}. 1.  Notice that $\Max (D)= \OO (h)\coprod \Max (T^{-1}D) $ where the inclusion $\Max (T^{-1}D) \subset \Max (D)$ is due to the injection $\Max (T^{-1}D)\ra \Max (D)$, $\gm\mapsto D\cap \gm$. Recall that every simple $D$-torsion $\CD (A)$-module is a simple weight $\CD (A)$-module, and vice versa, see (\ref{AD=AW}). Now, statement 1 follows from 
Lemma \ref{a26Jan20} and Proposition \ref{A29Jan20}. 

2. Statement 2 follows from statement 1 and Proposition \ref{A29Jan20}.

3. Statement 3 follows from statement 2. $\Box $\\

In order to describe the set of simple $D$-torsion free $\CD (A)$-modules we need to know a classification of  simple weight $\mA$-modules (Theorem \ref{B29Jan20})  and how  simple weight $\CD (A)$-modules with support from the orbit $\OO (h)$ decompose under restriction to the subalgebra $\mA$ of $\CD (A)$ (Lemma \ref{a29Jan20}). \\

\textbf{The set $\widehat{\mA}\,(D\text{-}{\rm torsion})=\widehat{\mA}\,({\rm weight})$}. Recall that the algebra $\mA$ is a generalized Weyl algebra $\mA = D[ \d_1, \d_{-1}; \s , a= h(h-1)(h-m)]$  where $D=K[h]$ and $\s (h) = h-1$. The orbit $\OO (h)$ is the only degenerated orbit and the maximal  ideals $(h)<(h-1)<(h-m)$ are the only marked maximal ideals. They partition the orbit $\OO (h)$ into subsets (see (\ref{GGG})):
$$\G_-=(-\infty, (h)], \,\,\,\,  \G_1=((h), (h-1)],\; \; \G_{m-1}=((h-1), (h-m)], \,\,\,\,\G_+=((h-m), \infty). $$

\begin{theorem}\label{B29Jan20} 
{\sc (Classification of simple $D$-torsion/weight $\mA$-modules)}  \\
The map $\Max\,(D)/\!\!\!\sim  \,\, \rightarrow \widehat{\mA}\,(D\text{-}{\rm torsion})$, $\G \mapsto [L(\G)]$, is  a bijection with the inverse $[M] \mapsto \Supp\,(M)$ where
\begin{enumerate}
\item If $\G\in \Max (D)/\langle \s \rangle\backslash \{ \OO (h)\}$ is  a non-degenerated orbit then $L(\G)=\mA/\mA\gp$ where $\gp \in \G$.
\item If $\G=\G_-=(-\infty, (h)]$ then $L_-:=L(\G_-)=\mA/\mA( h, \d_1)$.
\item If $\G=\G_1, \G_{m-1}$ then $L_1:=L(\G_1)=\mA/\mA(\d_{-1}, h-1 ,\d_1)$ and $L_{m-1}:=L(\G_{m-1})=\mA/\mA(\d_{-1}^{m-1}, h-m ,\d_1)$. These two modules are the only finite dimensional simple $\mA$-modules;  $\dim_K\, L(\G_1)=1$ and   $\dim_K\, L(\G_{m-1})=m-1$. 
\item If $\G=\G_+$ then $L_+:=L(\G_+)=\mA/\mA(h-m-1, \d_{-1})$.
\end{enumerate}
\end{theorem}
 
{\it Proof.} This is a particular case of Theorem \ref{AA27Feb16}. $\Box$\\

 Recall that $\mA \subset \CD (A) \subset \mB$. So every $\mB$-module is also an $\mA$-module and a $\CD (A)$-module (by restriction). Corollary \ref{b29Jan20} shows that the algebras  $\mA$, $ \CD (A)$ and $ \mB$ have the same simple $D$-torsion modules provided their supports do not belong to the orbit $\OO (h)$. For the algebras $R=\mA$, $\CD (A)$, $ \mB$, we denote by $\widehat{R } (D-{\rm torsion}\, | \, \OO )$ the set of simple $D$-torsion $R$-modules with support disjoint from $\OO (h)$.

  \begin{corollary}\label{b29Jan20}
$\widehat{\mA } (D-{\rm torsion}\, | \, \OO )=\widehat{\CD (A)} (D-{\rm torsion}\, | \, \OO )=\widehat{\mB } (D-{\rm torsion}\, | \, \OO )=\{ \mB / \mB \gp_\OO  \, | \, \gp_\OO  \in \Max (D)\backslash \OO (h)\}$ and  $\Supp (\mB / \mB \gp_\OO )= \OO $. 
\end{corollary}

{\it Proof}.  The corollary follows from the classifications of simple $D$-torsion modules for the algebras $ \mA$, $\CD (A)$ and $ \mB$ (Theorem \ref{DTor27Jan20} and  Theorem \ref{B29Jan20}). $\Box$

By Lemma   \ref{a26Jan20}, the $\CD (A)$-modules $A$ and $A'$ are the only two (up to isomorphism) simple weight $\CD (A)$-modules with support in the orbit $\OO (h)$. Lemma \ref{a29Jan20} shows that these modules are semisimple $\mA$-modules of length 2. 
 
 \begin{lemma}\label{a29Jan20}
\begin{enumerate}
\item ${}_{\mA} A= L_1\oplus L_+$ is a direct sum of simple weight $\mA$-modules  where the $\mA$-modules  $L_1$ and $ L_+$ are defined in Theorem \ref{B29Jan20}.(3,4). 
\item ${}_{\mA} A'=  L_-\oplus L_{m-1}$   is a direct sum of simple weight $\mA$-modules where the $\mA$-modules  $L_1$ and $ L_+$ are defined in Theorem \ref{B29Jan20}.(2,3). 
\end{enumerate}
\end{lemma}

{\it Proof}. 1. Recall that ${}_{\CD (A)} A=K+\sum_{i\geq m} Kx^i$. Then ${}_\mA K\simeq L_1$ and ${}_\mA \Big( \sum_{i\geq m} Kx^i\Big)\simeq L_+$.  Hence, ${}_{\mA} A= L_1\oplus L_+$ since $\d_1*K=0$, $\d_{-1}*K=0$ and $\d_{-1}*x^m=0$. 

2. Similarly,  ${}_{\CD (A)} A'=\Big(\sum_{i\leq -1} Kx^i\Big)\oplus \Big(\sum_{1\leq i\leq m-1} Kx^i\Big)$. Then $_\mA \Big(\sum_{i\leq -1} Kx^i\Big)\simeq L_-$ and ${}_\mA \Big(\sum_{1\leq i\leq m-1} Kx^i\Big)\simeq L_{m-1}$. Hence,  ${}_{\mA} A'= L_{m-1}\oplus L_-$  since $\d_1*x^{-1}=0$, $\d_{-1}*x=0$ and $\d_1*x^{m-1}=0$.  
 $\Box $ \\

 {\bf Classification of simple $D$-torsion free  $\CD (A)$-modules.}
 Recall that the algebra $\mA$ is a GWA $\mA = D[\d_1, \d_{-1}; \s , a=h(h-1)(h-m)]$. In order to stress that we consider `normal'  elements for the GWA $\mA$ we say `$\mA$-{\em normal}', see Theorem \ref{DTFr27Jan20}, i.e.   an element $b=\d_{-1}^m\beta_{-m}+\d_{-1}^{m-1}\beta_{-m+1}+\cdots + \beta_0 \in \mA$ (where $m \geqslant 1$, all $\beta_i \in D$ and $\beta_{-m}, \beta_0 \neq 0$) is called an $\mA$-{\em normal} element if $\beta_0 < \beta_{-m}$ and $\beta_0 < a$.

 \begin{theorem}\label{DTFr27Jan20}
{\sc (Classification of simple $D$-torsion free $\CD (A)$-modules)} 
 
$\widehat{\CD (A)}\,(D\text{-}{\rm torsion\; free})=\{ [M_b:=\CD (A)/\CD ( A)\cap Bb]\,|\,b$ is an $\mA$-normal irreducible element of $B\}$. The $\CD (A)$-modules $M_b$ and $M_{b'}$ are isomorphic iff the elements $b$ and $b'$ are similar.
\end{theorem}

{\it Proof}. Let $\CR$ be the RHS of the equality in the theorem.

(i) $\CR \subseteq \widehat{\CD (A)}\,(D\text{-}{\rm torsion\; free})$:
 Given $M_b:=[\CD (A)/\CD ( A)\cap Bb]\in \CR$ where $b$ is an $\mA$-normal irreducible element of $B$. We have to prove that $M_b\in \widehat{\CD (A)}\,(D\text{-}{\rm torsion\; free})$. By the very definition the $\CD (A)$-module $M_b$ is $D$-torsion free (since $M_b\subseteq B/ Bb$). By Theorem \ref{C29Jan20}, the $\CD (A)$-module $M_b$ has finite length. It remains to show that the $\CD (A)$-module $M_b$ is simple. Suppose that this is not the case, i.e. the left ideal $\CD (A) \cap Bb$ of the algebra $\CD (A)$  is not a maximal left ideal, we seek a contradiction. Then there is an element $\alpha \in D\backslash K$ such that the left ideal  $\CD (A) \cap Bb$ is properly contained in the left ideal $D\alpha +\CD (A) \cap Bb \neq \CD (A)$. Hence, let $W$ be a simple weight $\CD (A)$-factor module of the $\CD (A)$-module $\CD (A)/ (D\alpha +\CD (A) \cap Bb )$. In particular the action of the element $b\in \mA \subseteq \CD (A)$ has nonzero kernel. By Corollary \ref{b29Jan20} and   Lemma \ref{a29Jan20}, the weight  $\mA$-module $W$ is either simple or a direct sum of two  simple weight $\mA$-modules. 
 Hence, the action of the element $b$ is not injective on a simple $\mA$-submodule of $W$, this contradicts to  \cite[Lemma 3.7]{Bav-GWArep} since the element $b$ is $\mA$-normal, a contradiction.

(ii) $\CR \supseteq \widehat{\CD (A)}\,(D\text{-}{\rm torsion\; free})$:
 Let $M\in \widehat{\CD (A)}\,(D\text{-}{\rm torsion\; free})$. We have to show that $M\simeq M_b$ for some $\mA$-normal  irreducible element $b$  of $B$. The $B$-module $D^{-1}M$ is simple. Hence, $M\simeq M_b$ for some irreducible element of the algebra $B$. Since $D^{-1}\mA = B$ we may assume that $b=\d_{-1}^m\beta_{-m}+\d_{-1}^{m-1}\beta_{-m+1}+\cdots + \beta_0$ with all $\beta_i\in D$,  $\beta_{-m}\neq 0$ and $\beta_0\neq 0$.   By Lemma \ref{a27Feb16},   we may assume that the element $b$ is $\mA$-normal since the $B$-modules $B/Bb$ and $ B/B\beta b\alpha =  B/Bb \alpha$ are isomorphic (via the map $u\mapsto u\alpha $), and the statement (ii) follows. $\Box $


\section{The algebras  $\CDAm$ where $m\in \N^n$}\label{CDAMN}

In this section, properties of the algebras $\CDAm$ of differential operators  are studied where $m\in \N^n$. Proofs of Theorems \ref{GenProp30Jan20}--\ref{Kdim30JAn20} are given. The key idea of the proofs is to use properties of the generalized Weyl algebras $\CA (m)$ of rank $n$.

{\bf Generalized Weyl algebras of rank $n$, \cite{Bav-GWA-FA-91}-\cite{Bav-Bav-Ideals-II-93}.} Let $D$ be a ring, $\sigma=(\sigma_1,\ldots  ,\sigma_n)$   an $n$-tuple of
commuting automorphisms of $D$,  $a=(a_1,\ldots  ,a_n)$   an $n$-tuple  of elements of  the centre
$Z(D)$  of $D$ such that $\sigma_i(a_j)=a_j$ for all $i\neq j$. The {\bf  generalized Weyl algebra} $A=D(\sigma,a)=D[x, y; \sigma,a]$ of rank  $n$  is  a  ring  generated
by $D$  and    $2n$ indeterminates $x_1, \ldots ,x_n, y_1,\ldots  ,y_n$
subject to the defining relations:
$$y_ix_i=a_i,\;\; x_iy_i=\sigma_i(a_i),\;\; x_id=\sigma_i(d)x_i,\;\;{\rm and}\;\;  y_id=\sigma_i^{-1}(d)y_i\;\; {\rm for \; all}\;\; d \in D,$$
$$[x_i,x_j]=[x_i,y_j]=[y_i,y_j]=0, \;\; {\rm for \; all}\;\; i\neq j,$$
where $[x, y]=xy-yx$. We say that  $a$  and $\sigma $ are the  sets  of
{\it defining } elements and automorphisms of the GWA $A$, respectively.

The GWA $A=\bigoplus_{\alpha\in \Z^n}A_\alpha$ is a $\Z^n$-graded algebra ($A_\alpha A_\beta \subseteq A_{\alpha +\beta}$ for all elements $\alpha , \beta \in \Z^n$) where $A_\alpha = Dv_\alpha = v_\alpha D$, $v_\alpha =v_{\alpha_1}(1) \t\cdots \t v_{\alpha_n}(n)$, $v_{m}(i):=x_i^m$ and  $v_{-m}(i):=y_i^m$ for  all $m\geq 1$, and $v_0(i):=1$. \\

{\bf Example.} Let $D_i[x_i, y_i; \s_i , a_i]$ be GWAs of rank 1 over a field $K$ where $i=1, \ldots , n$. Then their tensor product over the field $K$,
$$\bigotimes_{i=1}^nD_i[x_i, y_i; \s_i , a_i]=D[x,y; \s , a],$$
is a GWA of rank $n$ where the $D=\bigotimes_{i=1}^n D_i$, $\s = (\s_1, \ldots , \s_n)$ and $a= (a_1, \ldots , a_n)$. The $\Z^n$-grading of the GWA $D[x,y; \s , a]$ of rank $n$ is the tensor product of $\Z$-gradings of the tensor components/GWAs of rank 1.

{\bf Example.}  The $n$'th Weyl algebra 
 $A_n=A_n(K)$ is a generalized Weyl algebra
 $A=D_n[x, y; \s ;a]$ of rank $n$ where
$D_n=K[h_1,...,h_n]$ is a polynomial algebra   in $n$ variables with
 coefficients in $K$, $\s = (\s_1, \ldots , \s_n)$ where $\s_i(h_j)=h_j-\delta_{ij}$,   $\d_{ij}$ is the Kronecker delta function and
 $a=(h_1, \ldots , h_n)$.  The map
$$A_n\ra A,\;\; x_i\mapsto  x_i,\;\; \der_i \mapsto y_i   ,\;\;  i=1,\ldots ,n,$$
is an algebra  isomorphism (notice that $\der_ix_i\mapsto h_i$). In particular, the GWA $A_n=\bigoplus_{\alpha \in \Z^n}D_nv_\alpha$ is a $\Z^n$-graded algebra where $v_\alpha =v_{\alpha_1}(1) \t\cdots \t v_{\alpha_n}(n)$, $v_{m}(i):=x_i^m$ and  $v_{-m}(i):=\der_i^m$ for all $m\geq 1$, and $v_0(i):=1$. \\

{\bf Generators and defining relations for the algebra $\CDAm$.} {\bf Proof of Theorem \ref{GenRel30Jan20}.} 

 1. The set $S_{n,x} :=\{ \prod_{i=1}^nx_i^{n_i} \; | \;  n_i \geq 0 \}$ 
(resp., $S_{n,x^m} := \{\prod_{i=1}^nx_i^{m_in_i} \; | \;   n_i \geq 0 \}$) 
is a multiplicative set of the polynomial algebra $P_n=K[x_1, \ldots , x_n] $ 
(resp., $P_n$ and $\Am$).
Clearly, 
\begin{equation}\label{SxmPAn1}
K[x,x^{-1}]:=K[x_1^{\pm 1}, \ldots , x_n^{\pm n}] = S_{n,x}^{-1} P_n = S_{n,x^m}^{-1} P_n = S_{n,x^m}^{-1} \Am .
\end{equation}
The set $S_{n,x}$ 
(resp., $S_{n,x^m}$) 
is an Ore set of the Weyl algebra $A_n$ 
(resp., of $A_n$ and $\CDAm$)
and 
\begin{equation}\label{ASADn1}
A_{n,x} := S_{n,x}^{-1} A_n = S_{n,x^m}^{-1} A_n = S_{n,x^m}^{-1} \CD(P_n) \simeq  \CD (S_{n,x^m}^{-1} P_n) \stackrel{(\ref{SxmPAn1})} = \CD (S_{n,x^m}^{-1} \Am ) \simeq  S_{n,x^m}^{-1} \CDAm .
\end{equation}
Recall that the Weyl algebra 
$A_n = D_n [x, \der ; \s , a  ]$ is a GWA of rank $n$, see above.  
In particular, the Weyl algebra 
$A_n = \bigoplus_{\alpha \in \Z^n} D_n v_\alpha$ is a $\Z^n$-graded algebra. 

Since the elements of the Ore set $S_{n,x}$ 
are homogeneous elements of the algebra $A_n$, 
the localized algebra $A_{n,x}$ is also a $\Z^n$-graded algebra
\begin{equation}\label{Anx}
A_{n,x} = \bigoplus_{\alpha \in \Z^n} D_n x^\alpha =D_n [x_1^{\pm 1},  \ldots ,  x_n^{\pm 1}; \s_1, \ldots , \s_n]
\end{equation}
which is a skew Laurent polynomial algebra where $D_n=K[h_1, \ldots , h_n]$, $h_i=\der_i x_i$ and  $x^\alpha = \prod_{i=1}^n x_i^{\alpha_i}$. 
By (\ref{ASADn1}),
$\CDAm = \{ \d \in A_{n,x} \; | \; \d \ast \Am \subseteq \Am  \}$.

Since the algebra $\Am$ is a $\Z^n$-graded subalgebra of the polynomial algebra $P_n$, the algebra $\CDAm$ is also $\Z^n$-graded,
\begin{equation}\label{ASAD1n1}
\CDAm = \bigoplus_{\alpha \in \Z^n} \CDAm_{[\alpha]}  \;\; {\rm  where} \;\;  \CDAm_{[\alpha]} = \CDAm \cap D_n x^\alpha  = \{ \d \in D_n x^\alpha  \; | \;  \d \ast \Am \subseteq \Am \}.
\end{equation}

Now, using the fact that $ h_i \ast x^{\alpha} = (\alpha_i+1) x^{\alpha}$ and  
for all $i=1, \ldots, n$ and  $\alpha \in \Z^n$, 
we obtain the explicit expressions for the graded components,   
$\CDAm_{[\alpha]}=\bigotimes_{i=1}^n \CD (A(m_i))_{[\alpha_i]}$, i.e. 
$\CDAm =\bigotimes_{i=1}^n \CD (A(m_i))$. 

2. Statement 2 follows from statement 1 and Theorem \ref{AB17Nov18}.(4).  $\Box$ \\

{\bf The algebra $\CDAm$ is a $\Z^n$-graded algebra.}  Recall that $\CDAm =\bigotimes_{i=1}^n\CD (A (m_i)) $ (Theorem \ref{GenRel30Jan20}.(1)).
If $m_i\geq 2$ then the algebra $\CD (A (m_i)) =\bigoplus_{j\in \Z}K[h_i]\d_j(i)$ is a $\Z$-graded algebra where the elements $\d_j(i)=\d_j$ are defined in Theorem \ref{AB18Nov18}.(1). If $m_i=1$ then the algebra $\CD (A(1))$ is the Weyl algebra $A_1=\bigoplus_{j\in \Z} K[h_i]\d_j(i)$ which is a $\Z$-graded algebra since it is a GWA where $\d_j(i)=\d_1(i)^j=x_i^j$ and $\d_{-j}(i)=\d_{-1}(i)^j=\der_i^j$ for $j\geq 0$. Since $\CDAm =\bigotimes_{i=1}^n\CD (A (m_i)) $ and every tensor component is a $\Z$-graded algebra the algebra $\CDAm$ is a $\Z^n$-graded algebra
\begin{equation}\label{grCDAm}
\CDAm =\bigotimes_{\alpha \in \Z^n} D_n\d_\alpha , \;\; D_n=K[h_1, \ldots , h_n], \;\; \d_\alpha =\prod_{i=1}^n\d_{\alpha_i}(i). 
\end{equation}
Notice that $D_n\d_\alpha = \d_\alpha D_n$ since $\d_\alpha d= \s^\alpha (d)\d_\alpha$ where $\s^\alpha = \prod_{i=1}^n \s_i$, $\s_i(h_j) = h_j-\d_{ij}$. The $\Z^n$-grading on the algebra $\CDAm$ in (\ref{grCDAm}) coincides with the induced $\Z^n$-grading that is determined by the embedding $\CDAm\subseteq A_{n,x}$ and the $\Z^n$-grading of the algebra $A_{n,x}$ in (\ref{Anx}).  \\

 {\bf The generalized Weyl algebras $\mA_n$ and $\mB_n$ such that $\mA_n \subset \CDAm \subset \mB_n \subset T_n^{-1}\mA_n =T_n^{-1} \CDAm = T_n^{-1} \mB_n $.} Recall that $\CDAm = \bigotimes_{i=1}^n \CD (A(m_i))$. For each number $i=1, \ldots , n$,   let $\mA (i)$ be the subalgebra of $\CD (A(m_i))$ which is generated by the elements $\d_{-1}(i)$, $h_i$ and $\d_1(i)$. By Theorem \ref{AB17Nov18}.(1), 
 $\d_{-1}(i)= h_i(h_i-m)x_i^{-1}$ and $\d_1(i)= (h_i-2)x_i$, and so the algebra 
\begin{equation}\label{ADddn1}
\mA (i) =D(i)[\d_1(i), \d_{-1}(i); \s_i , h_i(h_i-1)(h_i-m_i)], \;\; D(i)=K[h_i], \;\; \s_i (h_i) = h_i-1, 
\end{equation}
 is a GWA such that $\mA (i) \subset \CA_1(i) := \CD (A(m_i))\cap A_1(i)$ where $ A_1(i)=K\langle x_i , \der_i \, | \, \der_i x_i-x_i\der_i=1\rangle$ is the (first) Weyl algebra since $\d_{-1}(i), h_i , \d_1(i)\in  \CA_1(i)$ (Theorem \ref{BB18Nov18}.(2)). Let $$\mA_n:=\bigotimes_{i=1}^n \mA (i)\;\; {\rm and} \;\; \CA_n:=\bigotimes_{i=1}^n \CA_1 (i).$$
 Then $\mA_n \subseteq \CA_n$.

  The multiplicative submonoid $T(i)= \langle h_i-j\, | \, j\in \Z\rangle$ of $D(i)$ is a (left and right) denominator set of the algebras $\mA (i)$, $
 \CA_1(i)$, $\CD (A(m_i))$ and $A_1(i)$ such that 
\begin{equation}\label{BBTn1}
T(i)^{-1}\mA (i)\simeq T(i)^{-1}\CA_1 (i) \simeq T(i)^{-1} \CD (A(m_i)) \simeq  T(i)^{-1}A_1(i)=:  \mB (i)= T(i)^{-1}D(i)[x_i,x_i^{-1}; \s_i ]
\end{equation}
where $T(i)^{-1}D(i)=K[h_i, (h_i-j)^{-1}]_{j\in \Z}$ and $  \s_i (h_i) = h_i-1 $. Let 
$$\mB_n := \bigotimes_{i=1}^n\mB (i)\;\; {\rm  and }\;\;\CA (m):=\bigotimes_{i=1}^n \CA (m_i)$$
 where $\CA (m_i)$ is a subalgebra of  $\CD (A (m_i))$ which is generated by the elements 
$h_i$, $X_i:= x_i^{m_i}$ and $Y_i:= \d_{-m_i}(i)$. The algebra $\CA (m_i)$  is a  GWA of rank 1, 
 $$\CA (m_i)= K[h_i][X_i,Y_i;\s_i^{m_i}, a_i=(h_i+m_i-1) \cdot (h_i-2) (h_i-3) \cdots (h_i-m_i)], $$ which  is a central simple  Noetherian domain
where $\s_i (h_i) = h_i-1$, see Theorem \ref{AB18Nov18}.(2).

We have the following diagram of algebras where the vertical lines denote containments of the algebras where $T_n:= T(1) \cdots T(n)$ is a denominator set of the corresponding algebras:
\begin{figure}[H]
\centering
\begin{tikzpicture}
\node (B)  at (0,2){$\mB_n= T^{-1}_n\mA_n = T^{-1}_n\CA_n  = T^{-1}_n \CD (A(m)) =  T^{-1}_nA_n = T^{-1}_nD_n[x_1^{\pm 1}, \ldots ,x_n^{\pm 1}; \s_1, \ldots, \s_n]$};
\node (C) at (-1,1) {$\CD(A(m))$};
\node (D)  at (1,1){$A_n$};
\node (E)   at (0, 0){ \qquad \qquad \qquad \qquad  $\CA_n=\CD (A(m) \cap A_n)$};
\node (F)   at (-1, 0){$\CA(m)$};
\node (G)   at (0,-1){$\mA_n$};
\draw [semithick, -] (B)--(C);
\draw [semithick, -] (B)--(D);
\draw [semithick, -] (C)--(F);
\draw [semithick, -] (E)--(G);
\draw [semithick, -] (C)--(E);
\draw [semithick, -] (D)--(E);
\end{tikzpicture} 
\caption{AADBB1} \label{AADBB1}
\end{figure}

\begin{proposition}\label{A10Jan20}
Let $m=(m_1, \ldots ,m_n)\in \N^n$. Then 
\begin{enumerate}
\item The subalgebra $\CA (m)$ of $\CDAm$ is a  GWA of rank $n$ which  is a central simple  Noetherian domain of Gelfand-Kirillov dimension $2n$. 
\item The algebra $\CDAm$ is a finitely generated  left and right $\CA (m)$-module,  
$$ \CDAm  = \sum_{ \{ \alpha \in \Z^n : |\alpha_1| < 2m_1, \ldots ,  |\alpha_n| < 2m_n \} } \CA (m) \d_\alpha  = \sum_{ \{ \alpha \in \Z^n : |\alpha_1| < 2m_1, \ldots ,  |\alpha_n| < 2m_n \} } \d_\alpha  \CA (m). $$
\end{enumerate}
\end{proposition}

{\it Proof}. 2. Statement 2 follows from the fact that $\CDAm =\otimes_{i=1}^n\CD (A(m_i))$ and Theorem \ref{AB18Nov18}.(3). 

1. By \cite[Proposition 1.3]{Bav-GWArep}, the GWA $\CA (m)$ a   Noetherian domain. By \cite[Theorem 4.5]{Bav-FilDimSimCrit-1996}, the GWA $\CA (m)$ a  simple algebra. The algebra $\CA (m)$ is central since the algebra $A_{n,x}$ is so and 
$$\CA (m) \subset S_{n,x^m}^{-1}\CA (m) \simeq D_n[x_1^{\pm m_i}, \ldots , x_n^{\pm m_n}; \s_1^{m_1}, \ldots , \s_n^{m_n}]\simeq A_{n,x}\;\; (x_i^{m_i}\mapsto x_i, \;\; h_i\mapsto m_i h_i).$$
 The GWA $\CA (m)=\otimes_{i=1}^n\CA (m_i)$ is a tensor product of simple GWAs (see Theorem \ref{AB18Nov18}.(2)). By \cite[Corollary 4.8.(2)]{Bav-FilDimSimCrit-1996}, the Gelfand-Kirillov dimension of the algebra $\CA (m)$ is $2n$.  Now, by statement 2 and \cite[Propositiion 8.2.9.(ii)]{MR}, $\GK (\CDAm ) = \GK (\CA (m))=2n$. $\Box$ \\

{\bf Proof of Theorem \ref{GenProp30Jan20}.} 

 (i) {\em The algebra $\CDAm$ is central}: $K\subseteq Z(\CDAm )\stackrel{ (\ref{ASAD1n1})}{\subseteq} Z(A_{n,x})=K$, and so the  algebra $\CDAm$ is central. 

(ii) {\em The algebra $\CDAm$ is Noetherian with Gelfand-Kirillov dimension $2n$}: By Proposition \ref{A10Jan20}, the subalgebra $\CA (m)$ of $\CDAm$  is a Noetherian algebra of  Gelfand-Kirillov dimension $2n$ such that the algebra $\CDAm$ is a finitely generated left and right $\CA (m)$-module. Hence, the algebra $\CDAm$ is also Noetherian and by \cite[Proposition 8.2.9.(ii)]{MR},  $\GK (\CDAm ) = \GK (\CA (m))=2n$.  

(iii) {\em The algebra $\CDAm$ is simple and $\Z^n$-graded}: By (\ref{grCDAm}), the algebra $\CDAm$ is a $\Z^n$-graded algebra and the $\Z^n$-graded components $D_n\d_\alpha $ ($\alpha \in \Z^n$) of the algebra $\CDAm$ are the common eigenspaces of the commuting inner derivations $\ad_{h_1}, \ldots , \ad_{h_n}$ of the algebra $\CDAm$. Therefore every nonzero ideal of the algebra $\CDAm$ is a homogeneous ideal and as a result  has nontrivial intersection with the subalgebra $D_n$ of $\CDAm$. Since $D_n\subseteq \CA (m)$ and the algebra $\CA (m)$ is simple (Proposition \ref{A10Jan20}.(1)), all nonzero ideals of the algebra $\CDAm$ are equal to $\CDAm$, and so the algebra $\CDAm$ is a simple algebra.  $\Box$ \\

{\bf The Krull dimension of the algebras $\CD (\Am )$.} {\bf Proof of Theorem \ref{Kdim30JAn20}.} 
 
 By \cite[Corollary 4.8.(5)]{Bav-FilDimSimCrit-1996}, the Krull dimension of the GWA $\CA (m)$ is $n$.  By Proposition \ref{A10Jan20}, the algebra $\CDAm$ is  a finitely generated left and right $\CA (m)$-module. Hence, the Krull dimension of the algebra $\CDAm$ is smaller or  equal to the Krull dimension of the algebra $\CA (m)$ which is $n$. The polynomial algebra $D_n$ is the zero graded component of the $\Z^n$-graded algebra $\CDAm$. Hence, them map $I\mapsto \CDAm\t_{D_n}I$  (resp., $I\mapsto I\t_{D_n} \CDAm$) from the set of ideals of the algebra $D_n$ to the set of left (resp., right) ideals of the algebra $\CDAm$ is an injection. Hence,  the Krull dimension of $D_n$, which is $n$,  is smaller or equal to the Krull dimension of $\CDAm$. Therefore, the Krull dimension of the algebra $\CDAm$ is $n$. $\Box$\\

{\bf An analogue of the Inequality of Bernstein for the algebras  $\CD (\Am )$.}  By \cite[Corollary 4.8.(4)]{Bav-FilDimSimCrit-1996}, an analogue of the  the Inequality of Bernstein holds for the algebra $\CA (m)$: {\em For  all nonzero finitely generated $\CA (m)$-modules} $M$, $\GK (M)\geq n$.\\

 {\bf Proof of Theorem \ref{BerIn30JAn20}.}  By Proposition \ref{A10Jan20}, the algebra $\CDAm$ is  a finitely generated left and right $\CA (m)$-module.  Hence, each finitely generated $\CDAm$-module $M$ is a also a nonzero finitely generated $\CA (m)$-module. Now,
 $$ \GK_{\CDAm}(M)\geq  \GK_{\CA (m)}(M)\geq n,$$
and Theorem \ref{BerIn30JAn20} follows.  $\Box$\\

{\bf The global dimension of the algebras $\CD (\Am )$.} Recall that Morita equivalent algebras have the same global dimension and  the global  dimension of the Weyl algebra $A_n$ is $n$ (in characteristic zero). \\

{\bf Proof of Theorem \ref{gldim30JAn20}.}  By \cite[Theorem, p.29]{Musson-1986}, in the case $n=1$, the algebra $\CD (A(m_1))$ is Morita equivalent to the Weyl algebra $A_1$. Hence, for an arbitrary $n\geq 1$, the algebra $$\CD (A(m))=\CD (\otimes_{i=1}^nA(m_i))\simeq \otimes_{i=1}^n\CD (A(m_i))\;\; {\rm (Theorem \ref{GenProp30Jan20}.(1))},$$  where $m\in \N^n$, is Morita equivalent to the Weyl algebra $A_n=A_1^{\otimes n}$. 
 Therefore, the global dimension of the algebra $\CD (A(m))$ is equal to the global dimension of the Weyl algebra $A_n$ which is $n$. $\Box$

{\bf Licence.} For the purpose of open access, the author has applied a Creative Commons Attribution (CC BY) licence to any Author Accepted Manuscript version arising from this submission.

{\bf Disclosure statement.} No potential conflict of interest was reported by the author.

{\bf Data availability statement.} Data sharing not applicable – no new data generated.



\small{

}

\begin{tabular}{l  l}

V. V.   Bavula  \quad \quad \quad \quad  \quad \quad \quad \quad \quad \quad \quad  & K. H. Hakami \\ 
Department of Pure Mathematics \quad \quad \quad \quad  \quad \quad \quad \quad \quad \quad \quad  & Department of Mathematics\\
University of Sheffield \quad \quad \quad \quad  \quad \quad \quad \quad \quad \quad \quad  & Faculty of Science \\
Hicks Building \quad \quad \quad \quad  \quad \quad \quad \quad \quad \quad \quad  &  Jazan University \\
Sheffield S3 7RH \quad \quad \quad \quad  \quad \quad \quad \quad \quad \quad \quad  & Jazan 45142 \\
UK \quad \quad \quad \quad  \quad \quad \quad \quad \quad \quad \quad  &  Saudi Arabia \\
email: v.bavula@sheffield.ac.uk \quad \quad \quad \quad  \quad \quad \quad \quad \quad \quad \quad  & email: khakami@jazanu.edu.sa\\
\end{tabular}

\end{document}